\documentclass[11pt]{article}

\usepackage{a4wide}
\usepackage{amssymb}
\usepackage{amsmath,amsthm}
\usepackage{url}
\usepackage{graphicx,epsfig}
\usepackage{color}
\graphicspath{{Figures/}}
\input{psfig.sty}

\def\comment#1{}

\def\NI{\noindent}
\def\ni{\noindent}

\def\sk{\smallskip}
\def\cG{\widetilde{G}}

\def\thetasym{\Theta_{k,\ell}^{\circlearrowleft}}

\def\term#1{{\em #1}\marginpar{\raggedright{\small\it #1}}}
\def\ITEMMACRO #1 ??? #2 ???{\par\vskip4pt\noindent%
\hangindent=#2em\setbox0\hbox{#1\kern4pt}%
\ifdim\wd0<\hangindent\setbox0\hbox to\hangindent{\hss#1\kern7pt}\fi%
\box0\ignorespaces}

\def\Item(#1){\ITEMMACRO {\rm (#1)} ??? 1.8 ???}
\def\FreeItem#1{\ITEMMACRO {#1} ??? 1.8 ???}

\let\Bitem=\bItem
\def\BrackItem[#1]{\ITEMMACRO [#1] ??? 1.8 ???}

\def\Note{\ITEMMACRO {\NI\bf Note.}  ??? 1.8 ??? \\ \bgroup\sf}
\def\EndNote{\par\egroup\medskip\ni}

\newtheorem{theorem}{Theorem}[section]
\newtheorem{lemma}{Lemma}[section]
\theoremstyle{definition}\newtheorem{definition}{Definition}[section]
\newtheorem{proposition}{Proposition}[section]

\def\Fact#1.{\par\sk{\NI\bf Fact~#1.}\ }
\def\Problem#1.{\par\sk{\NI\bf Problem~#1.}\ }
\def\Claim#1.{\medbreak\ni{\bf Claim~#1.}\ }
\def\Case#1.{\medbreak\ni{\bf Case~#1.}\ }
\def\SubCase#1.{\medbreak\ni{\bf Subcase~#1.}\ }

\def\Proof{\ni{\sl Proof.}\ }
\def\qed{\hfill\fbox{\hbox{}}\medskip}
\def\qedclaim{\hfill$\triangle$\smallskip}

\definecolor{RED}{rgb}{.84,0,0}
\definecolor{BLUE}{rgb}{0,0,.75}
\def\CR{\color{RED}}
\def\CB{\color{BLUE}}
\def\BL{\color{black}}
\def\ovl{\overline}

\def\RR{\hbox{\sf I\kern-1ptR}}
\def\NN{\hbox{\sf I\kern-1ptN}}
\def\ZZ{\hbox{\sf Z\kern-4ptZ}}

\def\Min{\hbox{\sf Min}}
\def\Max{\hbox{\sf Max}}

\def\la{\nwarrow}
\def\ra{\nearrow}
\def\lb{\swarrow}
\def\rb{\searrow}

\def\fp{\alpha}
\def\fpl{\alpha_{\scriptscriptstyle \lb}}
\def\fpr{\alpha_{\scriptscriptstyle \ra}}
\def\fpla{\alpha_{\scriptscriptstyle \la}}
\def\fprb{\alpha_{\scriptscriptstyle \rb}}
\def\fplX{\hat{\alpha}_{\scriptscriptstyle \lb}}
\def\fprX{\hat{\alpha}_{\scriptscriptstyle \ra}}
\def\fX{\hat{\alpha}}

\def\bp{\beta}
\def\bX{\hat{\beta}}

\def\tsl{$--$}

\begin{document}
% *************************************************************
%                   TITLE PAGE
% *************************************************************
%

\title{{\bf\huge Bijections for Baxter Families
and Related Objects
}\\[6mm]}

\author{
      \normalsize Stefan Felsner%%
\footnote{Partially supported by DFG grant FE-340/7-1}\\
      \small\sf Institut f\"ur Mathematik,\\[-1mm]
      \small\sf Technische Universit\"at Berlin.\\[-1mm]
      \small\sf {\tt felsner@math.tu-berlin.de}
\and
      \normalsize \'Eric Fusy\\
      \small\sf Laboratoire d'Informatique (LIX)\\[-1mm]
      \small\sf \'Ecole Polytechnique.\\[-1mm]
      \small\sf {\tt fusy@lix.polytechnique.fr}
\and
      \normalsize Marc Noy\\
      \small\sf Departament de Matem\`atica Aplicada II\\[-1mm]
      \small\sf Universitat Polit\`ecnica de Catalunya.\\[-1mm]
      \small\sf {\tt marc.noy@upc.edu}
\and
      \normalsize David Orden%%
\footnote{Research partially supported by grants MTM2005-08618-C02-02 and S-0505/DPI/0235-02.}\\
      \small\sf Departamento de Matem\'aticas\\[-1mm]
      \small\sf Universidad de Alcal\'a \\[-1mm]
      \small\sf {\tt david.orden@uah.es}
}
\date{}
\maketitle
%%
%%

% *************************************************************

\vskip-0mm
\begingroup\fontsize{10}{12}\rm
\centerline{\normalsize\bf Abstract}
\bigskip
%abstract modified by Eric0303
\noindent The Baxter number $B_n$ can be written as $B_n = \sum_0^n \Theta_{k,n-k-1}$ with
$$
\Theta_{k,\ell} = \frac{2}{(k+1)^2\, (k+2)}
       {k+\ell \choose k}{k+\ell+1 \choose k}{k+\ell+2 \choose k}.
$$ 
%$$
%\Theta_{k,\ell} = 2\;\frac{(k+\ell)!\; (k+\ell+1)!\; (k+\ell+2)!}%%
%{k!\;(k+1)!\; (k+2)! \; \ell! \; (\ell+1)! \; (\ell+2)!}.
%$$ 
These numbers have first appeared in the enumeration of so-called
Baxter permutations; $B_n$ is the number of Baxter permutations of
size $n$, and $\Theta_{k,\ell}$ is the number of Baxter permutations
with $k$ descents and $\ell$ rises.  With a series of bijections we
identify several families of combinatorial objects counted by the
numbers~$\Theta_{k,\ell}$. Appart from Baxter permutations, these
include plane bipolar orientations with $k+2$ vertices and $\ell+2$
faces, 2-orientations of planar quadrangulations with $k+2$ white and
$\ell+2$ black vertices, certain pairs of binary trees with $k+1$ left
and $\ell+1$ right leaves, and a family of triples of non-intersecting
lattice paths. This last family allows us to determine the value of
$\Theta_{k,\ell}$ as an application of the lemma of Gessel and
Viennot. The approach also allows us to count certain other
subfamilies, e.g., alternating Baxter permutations, objects with
symmetries and, via a bijection with a class of plan bipolar orientations
also Schnyder woods of triangulations, which are known to be in bijection with
3-orientations.

\vskip10pt\ni
\endgroup

% *************************************************************
\section{Introduction}
% *************************************************************
\def\secline#1#2%%
{\medskip\hbox to\hsize{\bf\hbox to\parindent{\ref{sec:#2}\hfil}#1%%
\ \leaders\hbox{$\;.\;$}\hfill\ \pageref{sec:#2}}\smallskip}
\def\subsecline#1#2%%
{\smallskip\hbox to\hsize{\indent\hbox to 2em{{\bf \ref{ssec:#2}}\hfil}%%
#1\ \leaders\hbox{$\;.\;$}\hfill\ {\bf \pageref{ssec:#2}}}\smallskip}

This paper deals with combinatorial families enumerated by either
the Baxter numbers or the summands $\Theta_{k,\ell}$ of the usual
expression of Baxter numbers. Many of the enumeration results have
been known, even with bijective proofs. Our contribution to these
cases lies in the integration into a larger context and in
simplified bijections. We use specializations of the general
bijections to count certain subfamilies, e.g., alternating Baxter
permutations, objects with symmetries and Schnyder woods, i.e.,
3-orientations of triangulations.

This introduction will not include definitions of the objects we deal with, nor
bibliographic citations, which are gathered in notes throughout the article.
Therefore, we restrict it to a kind of commented table of contents.

{\advance\baselineskip-.5pt

\secline{Separating Decompositions and Book Embeddings}{dec_2book}

Separating decompositions of plane quadrangulations are defined. It is shown that
separating decompositions are in bijection with 2-orientations. Separating decompositions
induce book embeddings of the underlying quadrangulation on 2 pages. These special
book embeddings decompose into twin pairs of alternating trees, i.e.,
pairs of alternating trees with reverse reduced fingerprints. Actually, there is a bijection between
twin pairs of alternating trees and separating decompositions.

\secline{Alternating Trees and other Catalan Families}{altT+cat}

A bijection between alternating trees and full binary trees with the same fingerprint
is obtained. Fingerprint and bodyprint yield a bijection between full binary trees
with $k$ left and $\ell$ right leaves and certain pairs of non-intersecting lattice paths.
The lemma of Gessel and Viennot allows to identify their number as the Narayana number
$N(k+\ell-1,k)$.

\secline{Twin Pairs of Trees and the Baxter Numbers}{twins+baxt}

Twin pairs of alternating trees are in bijection with twin pairs of binary trees,
these in turn are shown to be in bijection to certain rectangulations and to
triples of non-intersecting lattice paths. Via the lemma of Gessel and Viennot
this implies that there are
$$
\Theta_{k,\ell} = \frac{2}{(k+1)^2\, (k+2)}
       {k+\ell \choose k}{k+\ell+1 \choose k}{k+\ell+2 \choose k}.
$$
twin pairs of binary trees with $k+1$ left and $\ell+1$ right leaves.
The bijections of previous sections yield a list of families enumerated
by the number $\Theta_{k,\ell}$.

\secline{More Baxter Families}{more_baxt}

We prove bijectively that
$\Theta_{k,\ell}$ counts Baxter permutations with $k$ descents and $\ell$ rises.
The bijections involve the Min- and Max-tree of a permutation and the
rectangulations from the previous section. Some remarks on the enumeration of
alternating Baxter permutations are added.

\subsecline{Plane Bipolar Orientations}{bipol}

We explain a bijection between separating decompositions and bipolar orientations.
The idea is to interpret the quadrangulation supporting the separating decompositions
as an angular map.

\subsecline{Digression: Duality, Completion Graph, and Hamiltonicity.}{hamiltonicity}

Combining ideas involving the angular map and the existence of a 2-book embedding for
quadrangulations we derive a Hamiltonicity result.

\secline{Symmetries}{symm}

The bijections between families counted by $\Theta_{k,\ell}$ have the nice property
that they commute with a half-turn rotation.  This is exploited to
count symmetric structures.

\secline{Schnyder Families}{schnyder}

Schnyder woods and 3-orientations of triangulations are known to be in bijection.
We add a bijection between Schnyder woods and bipolar orientations with a special property.
Tracing this special property through the bijections, we are able to
find the number of Schnyder woods on $n$ vertices via Gessel and Viennot.
This reproves a formula first obtained by Bonichon.

}%%% End of advance\baselineskip

% *************************************************************
\section{Separating Decompositions and Book Embeddings}\label{sec:dec_2book}
% *************************************************************

In the context of this paper a \term{quadrangulation}
is a plane graph $Q=(V\cup\{s,t\}, E)$ with only quadrangular faces.
More precisely, $Q$ is a maximal bipartite plane
graph, with $n+2$ vertices, prescribed color classes black and white and
two distinguished black vertices $s$ and~$t$ on the outer face.
Note that $Q$ has $n$ faces and $2n$ edges.

\begin{definition}
An orientation of the edges of $Q$ is a \term{2-orientation} if
every vertex, except the two special vertices $s$ and $t$, has
outdegree two. From easy arguments on double-counting the edges of $Q$,  $s$
and~$t$ are sinks in every 2-orientation.
\end{definition}

\begin{definition}\label{def:separating}
An orientation and coloring of the edges of $Q$ with colors red and
blue is a \term{separating decomposition} if:
\Item(1) All edges incident to $s$ are red and all edges incident to
$t$ are blue.
\Item(2) Every vertex $v\neq s,t$ is incident to an interval of red edges and
an interval of blue edges. If~$v$ is white, then, in clockwise order,
the first edge in the interval of a color is outgoing and all the
other edges of the interval are incoming. If~$v$ is white the outgoing
edge is the clockwise last in its color. (c.f.~Figure~\ref{fig:vertex-2-cond})
\end{definition}

%%%%%%%%%%%%%%%%%%%%%%%%%%%%%%%%%%%%%%%%%%%%%%%%%%%
%%
% in einem figure environment mit caption
   \calc_figscale{28}
    \begin{figure}[htb]
    \centerline{\input{\path/vertex-2-cond.pstex_t}}
    \caption{\label{fig:vertex-2-cond}}
    \end{figure}
    VC
{Edge orientations and colors at white and black vertices.}
%%
%%%%%%%%%%%%%%%%%%%%%%%%%%%%%%%%%%%%%%%%%%%%%%%%%%%

\begin{theorem}\label{thm:sd+2or}
Let $Q=(V\cup\{s,t\}, E)$ be a plane quadrangulation.
Separating decompositions and 2-orientations of $Q$ are in
bijection.
\end{theorem}

\Proof A separating decomposition clearly yields a 2-orientation,
just forget the coloring. For the converse, let $(v,w)$ be an
oriented edge, and define the \term{left-right path} of the edge
as the directed path starting with $(v,w)$ and taking a left-turn
in black vertices and a right-turn in white vertices.

\Claim A.
Every left-right path ends in one of the special vertices.
\smallskip

\ni{\it Proof of the claim.}
Suppose a left-right path closes a cycle $C$. The length of $C$ is an
even number $2k$. Let $r$ be the number of vertices interior to $C$.
Consider the inner quadrangulation of $C$. By Euler's formula it has
$2r+3k-2$ edges. However, when we sum up the outdegrees of the vertices we
find that $k$ vertices on $C$ contribute 1 while all other vertices
contribute 2 which gives a total of $2r + 3k$, contradiction.
\qedclaim

Now, let the special vertex where a left-right path ends
determine the color of all the edges along the path.

\Claim B. The two left-right paths starting at a vertex do not
meet again.
\smallskip

\ni{\it Proof of the claim.} Suppose that the two paths emanating
from $v$ meet again at $w$. The two paths form a cycle $C$ of even
length $2k$ with $r$ inner vertices. By Euler's formula the inner
quadrangulation of $C$ has $2r+3k-2$ edges. From the left-right
rule we know that one neighbor of $v$ on $C$ has an edge pointing
into the interior of $C$, from which it follows that there are at
least $k-1$ edges pointing from $C$ into its interior. Hence,
there are at least $2r + 3k- 1$ edges, contradiction. \qedclaim

Consequently, the two outgoing edges of a vertex~$v$ receive different colors.
It follows that the orientation and coloring of edges is a
separating decomposition.\qed

From the proof we obtain an additional property of a
separating decomposition:
\Item(3) \textit{The red edges form a directed tree rooted in $s$
         and the blue edges form a tree rooted in~$t$.}
\vskip4pt

\Note
De Fraysseix and Ossona de Mendez~\cite{fm-tao-01} defined
a separating decomposition via properties (1), (2) and (3), i.e.,
they included the tree-property into the definition.
They also proved Claim A and B and concluded
Theorem~\ref{thm:sd+2or}. In~\cite{fm-tao-01} it is also shown that
every quadrangulation admits a 2-orientation.
\EndNote

%\ni From the tree-property (3) it follows that, in a separating
%decomposition, every vertex~$v$ has a unique red path $P_r(v)$ to the red
%root~$s$ and a unique blue path $P_b(v)$ to the blue root~$t$. From the above
%claims we know that these two paths are simple and disjoint.  Let
%$P(v) = P_r(v)^{-1} \circ P_b(v)$ be the path from~$s$ to~$t$ through~$v$ which
%takes $P_r(v)$ in reverse orientation followed by $P_b(v)$.  Let
%$R(v)$ be the region left of $P(v)$, and define $v \prec w$ if $R(v)
%\subset R(w)$. Properties of separating trees imply that there is a
%containment of regions for any two vertices, i.e., the $\prec$-relation
%is a total order.

%\begin{proposition}\label{prop:2book}
%If the vertices of a quadrangulation $Q$ endowed with a separating decomposition
%are arranged on a line according to the $\prec$-ordering with
%$s$ being first and $t$ being last, then there is a plane
%drawing of $Q$ with all red edges below and all blue edges above the
%line.  In other words, the coloring of edges corresponds to the pages
%of a 2-book embedding of $Q$ with the referred line as spine. (c.f.~Figure~\ref{fig:sep2book})
%\end{proposition}

%\Proof
%The next lemma is an observation about the spine of a book embedding,
%which leads to our proof for the existence of a 2-book embedding of $Q$.

An embedding of a plane graph is called a \emph{2-book embedding}
if the vertices are arranged on a single line so that all edges are either 
below or above the line. As we show next, a separating decomposition $S$   
easily yields a 2-book embedding of the underlying quadrangulation $Q$.

Observe that each inner face of $Q$ has exactly two
bicolored angles, i.e., angles where edges of different color
meet; this follows from the rules given in Definition
\ref{def:separating}. Define the \term{equatorial line} of $S$  as the union of all diagonals connecting bicolored angles of
an inner face. The
definition of a separating decomposition implies that each inner
vertex of $Q$ has degree two in the equatorial line, while $s$ and
$t$ have degree zero and the other two outer vertices have degree
one. This implies that the equatorial line consists of an
edge-disjoint union of a path and possibly a collection of cycles.

\begin{lemma}\label{lem:once_faces}
Given a quadrangulation $Q$ endowed with a separating decomposition $S$, 
the equatorial line of $S$ consists of a single path that 
traverses every inner vertex and every inner face of $Q$ exactly once.
\end{lemma}

\Proof 
Assume that the equatorial line has a cycle $C$. Consider a
plane drawing of $Q\cup C$. The cycle $C$ splits the drawing into
an inner and an outer part, both special vertices $s$ and $t$
being in the outer part.  The red edges of all vertices of $C$
emanate to one side of $C$ while the blue edges go to the other
side. Therefore, it is impossible to have a monochromatic path
from a vertex $v\in C$ to both special vertices. With property~(3) of
separating decompositions, it thus follows that there are no
cycles, i.e., the equatorial line is a single path. 

 The equatorial line $L$ has  the two outer non-special vertices of $Q$ as extremities. 
It has degree 2 for each inner vertex, hence it traverses
 each inner vertex once. In addition, $L$  
separates the blue and the red edges, hence it must pass through the
interior of each inner face at least once. Since the $n$
non-special vertices of $Q$ delimit $n-1$ intervals on $L$ 
and since $Q$ has $n-1$ inner faces, $L$ can only pass
through each inner face exactly once. \qed

To produce a 2-book embedding, 
 extend the equatorial line in the outer face so that it visits
also the two special vertices; then stretch the equatorial line
as a straight horizontal line  that has $s$ as its leftmost vertex 
 and $t$ as its rightmost vertex, see Figure~\ref{fig:sep2bookWithEqLine}. 
%arrange the vertices of~$Q$ in the order in which
%they appear on the equatorial line, with $s$ as leftmost and
%$t$ as rightmost vertex on the line. 
Observe that 
one page gathers the blue edges, the other page gathers the red edges,
and the spine for the two pages is the equatorial line.

The equatorial line will be useful again for proving a Hamiltonicity
result in Section~\ref{ssec:bipol}.

%%%%%%%%%%%%%%%%%%%%%%%%%%%%%%%%%%%%%%%%%%%%%%%%%%%
%%
% in einem figure environment mit caption
   \calc_figscale{22}
    \begin{figure}[htb]
    \centerline{\input{\path/sep2bookWithEqLine.pstex_t}}
    \caption{\label{fig:sep2bookWithEqLine}}
    \end{figure}
    VC
{A quadrangulation $Q$ with a separating decomposition $S$, and the 
  2-book embedding induced by the equatorial line of $S$.}
%%
%%%%%%%%%%%%%%%%%%%%%%%%%%%%%%%%%%%%%%%%%%%%%%%%%%%
%%
\begin{definition}
An \term{alternating layout} of a tree $T$ with $n+1$ vertices is
a non-crossing drawing of $T$ such that its vertices are mapped to
different points on the non-negative $x$-axis and all edges are
embedded in the halfplane above the $x$-axis (or all below).
Moreover, for every vertex~$v$ it holds that all its neighbors are
on one side, either they are all left of~$v$ or all right of $v$.
In these cases we call the vertex~$v$ respectively a \term{right}
or a~\term{left vertex} of the alternating layout. A tree with an
alternating layout is an \emph{alternating tree}.
\end{definition}

\marginpar{\it alternating tree}

As one can check in the example of Figure~\ref{fig:sep2bookWithEqLine},
 the red and the blue trees
are both alternating. Indeed, by the rules given in Definition
\ref{def:separating}, white vertices are right in
the red tree and left in the blue tree while black vertices are left
in red and right in blue. 
%It has been shown in~\cite{fhko-blpqr-07}
%that this is always the case. 
We summarize:

\begin{proposition}\label{prop:alt-trees}
The 2-book embedding induced by a separating decomposition yields
simultaneous alternating layouts of the two trees such that each
white vertex is left in the blue tree and right in the red tree,
while each black vertex is right in the blue tree and left in the
red tree.
\end{proposition}

\Note
A proof of Proposition~\ref{prop:alt-trees} was given by
Felsner, Huemer, Kappes and Orden~\cite{fhko-blpqr-07}. These
authors study what they call {\sl strong binary labelings of the
angles of a quadrangulation}. They show that these labelings are
in bijection with 2-orientations and separating decompositions. In
this context they find the 2-book embedding; their method consists
in ranking each vertex $v$ on the spine of the 2-book embedding
according to the number of faces in a specific region $R(v)$. The original source
for a 2-book embedding of a quadrangulation
is~\cite{fmp-lrspg-95}, by de Fraysseix, Ossona de Mendez and
Pach. General planar graphs may require as many as 4 pages for a
book embedding, Yannakakis~\cite{y-epg4p-86}.

\FreeItem{}
Alternating trees in our sense were studied
by Rote, Streinu and Santos~\cite{rss-ecpmpt-03} as {\sl non-crossing alternating trees}.
There it is pointed out that the name {\em alternating tree} is
sometimes used to denote a tree with a numbering such that every
vertex is a local extremum, e.g., \cite[Exercise 5.41]{s-ec2-99}.
Non-crossing alternating trees were studied by Gelfand et
al.~\cite{ggp-chfapr-97} under the name of {\sl standard trees};
there it is shown that these trees are a Catalan
family. In~\cite{rss-ecpmpt-03} connections with rigidity theory and
the geometry of the associahedron are established.
\EndNote

\ni Consider alternating layouts of rooted ordered trees with the
property that the root is extreme, i.e., the root is either the
leftmost or the rightmost of the vertices. In this setting an
alternating layout is completely determined by the placement of
the root (left/right) and the choice of a halfplane for the edges
(above/below). We denote the four choices with symbols, e.g.,
$\lb$ denotes that the root is left and the halfplane below; the
symbols $\la$, $\ra$ and $\rb$ represent the other three
possibilities.

The unique $\lb$-alternating layout of $T$,
is obtained by starting at the root and walking
clockwise around $T$, thereby numbering the vertices with
consecutive integers according to the following rules: The root is
numbered $0$ and all vertices in the color class of the root receive a
number at the first visit while the vertices in the other color class
receive a number at the last visit. Figure~\ref{fig:alt-tree2} shows an
example. Rules for the other types of alternating layouts are:

\FreeItem{} $\la$-layout: walk counterclockwise, root class at first visit,
other at last visit.

\FreeItem{} $\ra$-layout: walk counterclockwise, root class at last
visit, other at first visit.

\FreeItem{} $\rb$-layout: walk clockwise, root class at last visit,
other at first visit.
\medskip

%%%%%%%%%%%%%%%%%%%%%%%%%%%%%%%%%%%%%%%%%%%%%%%%%%%
%%
% in einem figure environment mit caption
   \calc_figscale{22}
    \begin{figure}[htb]
    \centerline{\input{\path/alt-tree2.pstex_t}}
    \caption{\label{fig:alt-tree2}}
    \end{figure}
    VC
{A tree, the numbering and the $\lb$-alternating layout.}
%%
%%%%%%%%%%%%%%%%%%%%%%%%%%%%%%%%%%%%%%%%%%%%%%%%%%%

\ni The \term{$\lb$-fingerprint}, denoted $\fpl(T)$, of a rooted
ordered tree $T$, is a $0,1$ string which has a $1$ at
position~$i$, i.e., $\fp_i = 1$, if the $i$th vertex in the
$\lb$-alternating layout of $T$ is a left vertex, otherwise, if
the vertex is a right vertex $\fp_i = 0$. The $\lb$-fingerprint of
the tree $T$ from Figure~\ref{fig:alt-tree2} is $\fpl(T) = 1 0 1 0
0 0 1 0 1 0 0 0 0 1 1 0$.  Other types of fingerprints are defined
by the same rule.  For example the
%%$\la$-fingerprint of the tree in Figure~\ref{fig:alt-tree2} is
%%$\fp_{\scriptscriptstyle \la}(T) = 1 1 1 0 1 0 0 0 1 0 0 1 0 0 0 0$.
%%With the numbering from Figure~\ref{fig:alt-tree2} this corresponds to
%%the vertex order $0,14,13,15,8,11,10,9,6,7,12,2,4,3,5,1$.
$\ra$-fingerprint of the tree in Figure~\ref{fig:alt-tree2} is
$\fp_{\scriptscriptstyle \ra}(T) = 1 0 0 1 1 1 1 0 1 0 1 1 1 0 1 0$.
With the numbering from Figure~\ref{fig:alt-tree2} this corresponds to
the vertex order $0,15,14,13,12,11,10,9,8,7,6,5,4,3,2,1$.

The first vertex is always a left
vertex and the last a right vertex, therefore, a fingerprint
has always  a 1 as first entry and a 0 as last entry.
A \term{reduced fingerprint} $\fplX(T)$, of a tree $T$ is
obtained by omitting the first and the last entry from the
corresponding fingerprint.
For a $0,1$ string $s$ we define $\rho(s)$ to be
the reverse string and $\ovl{s}$ to be the complemented
string. Example: if $s=11010$, then $\rho(s) = 01011$, $\ovl{s}
= 00101$, and $\ovl{\rho(s)}=\rho(\ovl{s})=10100$.

\begin{lemma}\label{lem:fp-transformation}
For every tree $T$ we have $\fpl(T) = \ovl{\rho(\fpr(T))}$
(and $\fpla(T) = \ovl{\rho(\fprb(T))}$).
\end{lemma}

\Proof
Take the $\ra$-alternating layout of $T$ and rotate it by
$180^\circ$. This results in the $\lb$-alternating layout.
Observe what happens to the fingerprint.
\qed

\begin{definition}
A pair $(S,T)$ of rooted, oriented trees
whose fingerprints satisfy
$\fplX(S) = \ovl{\fprX(T)}$, or equivalently $\fprX(S) = \rho(\fprX(T))$, is called a
\term{twin-alternating pair of trees}.
\end{definition}

\begin{theorem}\label{thm:2or+ta}
There is a bijection between twin-alternating pairs of trees
$(S,T)$ on $n$ vertices and 2-orientations of quadrangulations
on $n+2$ vertices.
\end{theorem}

\Proof (The bijection is illustrated with an example in
Figure~\ref{fig:small-ex}.) Augment both rooted ordered trees $S$ and
$T$ by a new vertex which is made the rightmost child of the root. Let
$S^+$ and $T^+$ be the augmented trees. Note that $\fplX(S^+) = 0 +
\fplX(S)$ and $\fprX(T^+) = \fprX(T) + 1$. Since the first entry of a
non-reduced fingerprint is always 1 and the last one is always 0 it
follows that $\fpl(S^+) + 0 = \ovl{1 + \fpr(T^+)}$.

Consider the $\lb$-alternating layout of $S^+$ and move the vertices
in this layout to the integers $0,..,n$.  Similarly, the
$\ra$-alternating layout of $T^+$ is placed such that the vertices
correspond to the integers $1,..,n+1$. At every integer $0<i<n+1$ a
vertex of $S^+$ and a vertex of $T^+$ meet. We identify them. As a
consequence of the complemented fitting of the fingerprints, every
non-special vertex is a left vertex in one of the layouts and a right
vertex in the other. This has strong consequences:

\Bitem
A pair $uv$ can be an edge in at most one
of $S$ and $T$, otherwise $u$ would have a neighbor on its right in
both $S$ and $T$, a contradiction.

\Bitem There is no triangle with edges from $S\cup T$. Suppose
$u,v,w$ would be such a triangle. Two edges must be from the same
tree, say from $S$. These cannot be the two edges incident to the
middle vertex~$v$. If they are incident to $w$ the vertex $u$ has
neighbors to its right in both trees, contradiction.
\smallskip

\noindent
Hence, the graph with edges $S^+\cup T^+$ is simple, triangle-free and
non-crossing. Since it has $n+2$ vertices and
$2n$ edges, it must be a quadrangulation. The $2$-orientation is
obtained by orienting both trees towards the root.

%%%%%%%%%%%%%%%%%%%%%%%%%%%%%%%%%%%%%%%%%%%%%%%%%%%
%%
% in einem figure environment mit caption
   \calc_figscale{22}
    \begin{figure}[htb]
    \centerline{\input{\path/small-ex.pstex_t}}
    \caption{\label{fig:small-ex}}
    \end{figure}
    VC
{A twin-alternating pair of trees $(S,T)$.
 The $\ra$-alternating layout of $T^+$ and the $\lb$-alternating layout of $S^+$ properly adjusted.
 The induced 2-orientation of a quadrangulation.}
%%
%%%%%%%%%%%%%%%%%%%%%%%%%%%%%%%%%%%%%%%%%%%%%%%%%%%

The converse direction, from the 2-orientation of a quadrangulation on
$n+2$ vertices to trees $(S,T)$ with appropriate fingerprints was
already indicated.  To recapitulate: A 2-orientation of $Q$ yields a
separating decomposition (Theorem~\ref{thm:sd+2or}). In particular the
edges are decomposed into two trees, the red tree $S^+$ and the blue
tree $T^+$. The corresponding 2-book embedding
(Proposition~\ref{prop:alt-trees}) induces a simultaneous alternating
layout of the two trees with the property that every non-special
vertex is left in one of the trees and right in the other
(Proposition~\ref{prop:alt-trees}), i.e., $\fpl(S^+) + 0 = \ovl{1 +
  \fpr(T^+)}$.  Trees $S$ and $T$ are obtained by deleting the left
child of the root in $S^+$ and the right child of the root in $T^+$;
they are both leaves and correspond to the two non-special outer
vertices of $Q$. Trees $S$ and $T$ satisfy $\fplX(S) = \ovl{\fprX(T)}$,
i.e., $(S,T)$ is a twin-alternating pair of trees.  \qed

% *************************************************************
\section{Alternating Trees and other Catalan Families}\label{sec:altT+cat}
% *************************************************************

A \term{full binary tree} is a rooted ordered tree
such that each inner vertex has exactly two children.
The \term{fingerprint} of a full binary tree $T$ is
a $0,1$ string  which has a $1$ at position $i$ if the
$i$th leaf of $T$ is a left child,
otherwise, if the leaf is a right child the entry is $0$.
In Figure~\ref{fig:alt-bin} the tree~$T$ on the right side has
$\fp(T) = 1 0 1 1 1 0 1 0 1 1 1 1 0 0 1 0$.
The \term{reduced fingerprint} $\fX(T)$ is obtained by
omitting the first and the last entry from $\fp(T)$.
Note that the first entry is always~1 and the last one is always~0.

\def\lambdaT{T^\lambda}

\begin{proposition}\label{prop:fingerp-bij}
There is a bijection $T \to \lambdaT$ which takes an alternating
tree $T$ with $n$ vertices
to a full binary tree $\lambdaT$ with $n$ leaves
such that $\fprX(T) = \fX(\lambdaT)$.
\end{proposition}

\Proof
The bijection makes a correspondence between edges of the alternating
tree and inner vertices of the full binary tree, see
Figure~\ref{fig:alt-bin}. Embed $T$ with vertices on integers
from $0$ to~$n$. With an edge $i,j$ of $T$ associate an inner vertex $x_{ij}$
for $\lambdaT$ which is to be placed at
$(\frac{i+j}{2},\frac{j-i}{2})$. Draw line segments from the vertex
$(i,0)$ to $x_{ij}$ and from $(j,0)$ to $x_{ij}$. Doing this for every
edge of $T$ results in a drawing of the binary tree $\lambdaT$.
%%
%%%%%%%%%%%%%%%%%%%%%%%%%%%%%%%%%%%%%%%%%%%%%%%%%%%
%%
% in einem figure environment mit caption
   \calc_figscale{22}
    \begin{figure}[htb]
    \centerline{\input{\path/alt-bin.pstex_t}}
    \caption{\label{fig:alt-bin}}
    \end{figure}
    VC
{An $\ra$- alternating tree $T$ and the full binary tree $\lambdaT$.}
%%
%%%%%%%%%%%%%%%%%%%%%%%%%%%%%%%%%%%%%%%%%%%%%%%%%%%
%%

The converse is even simpler. Every inner node $x$ of the binary tree
gives rise to an edge connecting the leftmost leaf below $x$ to the rightmost
leaf below $x$.\qed

\Note
Full binary trees with $n+1$ leaves are counted by the
\term{Catalan number} $C_n = \frac{1}{n+1}\binom{2n}{n}$.
Catalan numbers are found in The On-Line Encyclopedia of
Integer Sequences~\cite{oleis} as
sequence~A000108.
From Proposition~\ref{prop:fingerp-bij} it follows that
alternating trees with $n+1$ vertices are another Catalan
family, which in~\cite{ggp-chfapr-97} was proved constructing a bijection inductively.
Stanley~\cite[Exercise 6.19]{s-ec2-99} collected 66 Catalan families.
\EndNote

\ni
Although the subject is well-studied, we include a particular
proof showing that full binary trees are a Catalan family.
Actually, we prove a more refined count related to Narayana numbers.
The proof will be used later in the context of Baxter numbers.

To start with, we associate another $0,1$ string with a full binary tree
$T$. The \term{bodyprint}~$\bp(T)$ of~$T$ is obtained from a
visit to the inner vertices of $T$ in in-order. The $i$th entry of $\bp$ is
a~$1$, i.e., $\bp_i = 1$, if the $i$th inner vertex is a right-child or it is the root.
If the vertex is a left-child, then $\bp_i=0$. Note that if the tree $T$
is drawn such that all leaves are on a horizontal line, then there is a
one-to-one correspondence between inner vertices and the gaps between adjacent leaves
(Gap between leaves $v_i$ and $v_{i+1}$ $\mapsto$ least common ancestor of $v_i$ and $v_{i+1}$.
Inner vertex $x$ $\mapsto$ gap between rightmost leaf in left subtree below~$x$ and
leftmost leaf in right subtree below~$x$). This correspondence maps the left-to-right
order of gaps between leaves to the in-order of inner vertices.
Since the root contributes a~$1$ the last entry of the bodyprint of a tree is always~$1$.
Therefore, it makes sense to define the \term{reduced bodyprint} $\bX(T)$
as $\bp(T)$ minus the last entry.
Figure~\ref{fig:enc-trees} shows an example.
%%
%%%%%%%%%%%%%%%%%%%%%%%%%%%%%%%%%%%%%%%%%%%%%%%%%%%
%%
% in einem figure environment mit caption
   \calc_figscale{82}
    \begin{figure}[htb]
    \centerline{\input{\path/enc-trees.pstex_t}}
    \caption{\label{fig:enc-trees}}
    \end{figure}
    VC
{A full binary tree with reduced bodyprint $\bX$ and reduced fingerprint $\fX$.}
%%
%%%%%%%%%%%%%%%%%%%%%%%%%%%%%%%%%%%%%%%%%%%%%%%%%%%
%%

\begin{lemma}\label{lem:a-dom-b}
Let $T$ be a full binary tree with  with $k$ left leaves and $n-k+1$ right leaves.
The reduced fingerprint $\fX(T)$ and the reduced bodyprint $\bX(T)$ both have length $n-1$.
Moreover:\\[2mm]
\centerline{
{\rm (1)}\quad $\displaystyle \sum_{i=1}^{n-1} \fX_i = \sum_{i=1}^{n-1} \bX_i = k-1$
\hskip14mm
{\rm (2)}\quad $\displaystyle \sum_{i=1}^{j} \fX_i \geq \sum_{i=1}^{j} \bX_i$\quad
         for all $j=1,\ldots,n-1$.
}
\end{lemma}

\Proof Consider a drawing of $T$ where every edge has slope $1$ or
$-1$, as in Figure~\ref{fig:enc-trees}. The maximal segments of slope
$1$ in this drawing define a matching $M$ between the $k$ left leaves, i.e.,
1-entries of $\fp(T)$, and inner vertices which are right-childs including
the root, i.e., 1-entries of $\bp(T)$. The left part of
Figure~\ref{fig:enc-trees-the-1s} indicates the correspondence.
The reduction $\fX$ (resp. $\bX$) has
exactly one 1-entry less than $\fp$ (resp. $\bp$). This proves (1).

%%%%%%%%%%%%%%%%%%%%%%%%%%%%%%%%%%%%%%%%%%%%%%%%%%%
%%
% in einem figure environment mit caption
   \calc_figscale{72}
    \begin{figure}[htb]
    \centerline{\input{\path/enc-trees-the-1s.pstex_t}}
    \caption{\label{fig:enc-trees-the-1s}}
    \end{figure}
    VC
{Illustrations for the proof of Lemma~\ref{lem:a-dom-b}.
 The pair $(v_i,x_j)$ is in $M$.}
%%
%%%%%%%%%%%%%%%%%%%%%%%%%%%%%%%%%%%%%%%%%%%%%%%%%%%

For (2) let $v_0,v_1,\ldots,v_n$ be the set of leaves in left-to-right
order and let $x_1,\ldots,x_n$ be the in-order of inner vertices. Note
that $v_i$ determines $\fp_i$ and $x_i$ determines $\bp_i$. Let
$(v_i,x_j)$ be a pair from the matching $M$ defined above,
i.e., $\fp_i = 1$ and $\bp_j = 1$. Since $v_i$ is the leftmost leaf
below $x_j$ and the gap corresponding to $x_j$ starts at the rightmost
vertex $v_{j-1}$ of the left subtree of $x_j$ we find that $i \leq
j-1$. This gives a matching between the 1-entries of $\fp$
and the 1-entries of $\bp$ with the property that the index of the
1-entry of $\fp$ is always less that the index of the mate in $\bp$.

To conclude the inequality for the reduced strings we have to address
another detail: The mate of the root in $M$ is the leaf $v_0$,
which is not represented in $\fX$, and there is a leaf whose mate
in $M$ is the last inner vertex $x_n$, which is not represented in
$\bX$.  Consider the ordered sequence
$x_{j_0},x_{j_1},\ldots,x_{j_s}$ of all vertices on the rightmost
branch of $T$, such that $x_{j_0}$ is the root $r$ and $x_{j_s} =
x_n$.  The right part of Figure~\ref{fig:enc-trees-the-1s} may
help to see that in $M$ we have the following pairs
$(v_0,x_{j_0}), (v_{j_0},x_{j_1}), \ldots (v_{j_{s-1}}, x_n)$; in
particular $\fp_{0} = \fp_{j_0} = \ldots =\fp_{j_{s-1}} = 1$ and
$\bp_{j_0} = \bp_{j_1} = \ldots =\bp_{n} = 1$. Hence we can define
a matching $M'$ which is as $M$ except that $v_0$ and $x_n$ remain
unmatched and the pairs $(v_{j_i},x_{j_i})$ with $0\leq i\leq s-1$
are matched.  This matching $M'$ between the 1-entries of $\fX$
and the 1-entries of $\bX$ has the property that the index of the
1-entry of $\fX$ is always at most the index of the mate in $\bX$.
This proves (2).  \qed

\def\geqdom{\geq_{{\sf dom}}}
\def\leqdom{\leq_{{\sf dom}}}
\begin{definition}
With $\sigma \in \binom{k+\ell}{k}$ we denote that $\sigma$ is a $0,1$ string of
length $n=k+\ell$ with $k$ entries $1$ and $\ell$ entries $0$, i.e.,
$\sum_{i=1}^{n} \sigma_i = k$.
For $\sigma,\tau \in \binom{k+\ell}{k}$ we define $\tau \geqdom \sigma$, i.e.,
$\tau$ \term{dominates} $\sigma$, if
$\sum_{i=1}^{j} \tau_i \geq \sum_{i=1}^{j} \sigma_i$
for all $j=1,\ldots, n$.
\end{definition}

\begin{theorem}\label{thm:tree-2strings}
The mapping $T \leftrightarrow (\bX,\fX)$ is a bijection between
full binary trees with $k+1$ left leaves and  $\ell+1$ right leaves and pairs $(\bX,\fX)$ of
$0,1$ strings in $\binom{k+\ell}{k}$ with $\fX \geqdom \bX$.
\end{theorem}

\Proof
From Lemma~\ref{lem:a-dom-b} we know that reduced body- and fingerprint have the
required properties. To show that the mapping $T \leftrightarrow (\bX,\fX)$ is a bijection
we use induction.

First note that $\fX = 0^{\ell}1^{k}$ implies $\bX=\fX$, and that there are
unique trees with these reduced finger- and bodyprints.

If $\fX$ has a different structure, then there is an $i$ such that
$\fX_{i-1}\fX_i = 10$. In a hypothetical tree $T$ corresponding
to $(\bX,\fX)$, the pair $v_{i-1},v_i$ is a left leaf
followed by a right leaf. The leaves $v_{i-1}$ and $v_i$ are children of
the inner vertex $x_i$, the value $\bX_i=1$ or $\bX_i=0$ depends on whether $x_i$ is
itself a left or a right child. Consider the tree $T^*$ obtained by pruning
the two leaves $v_{i-1}$ and $v_i$. Note that if $\fX(T) = \fX'\>\fX_{i-1}\>\fX_i\>\fX''$
and  $\bX(T) = \bX'\>\bX_i\>\bX''$, then $\bX(T^*) = \bX'\>\bX''$
and $\fX(T^*) = \fX'\>\delta\>\fX''$ where $\delta = 1$ if $\bX_i = 0$
and $\delta = 0$ if $\bX_i = 1$.

Hence, we define $\fX^* = \fX'\>\delta\>\fX''$ and $\bX^* =
\bX'\>\bX''$. Depending on the value of $\delta = \ovl{\bX_i}$, this can be interpreted
as either having removed the two entries $\bX_i=1$ and $\fX_{i-1}=1$
or the two entries $\bX_i=0$ and $\fX_{i}=0$ from $\fX$ and~$\bX$.  It
is easy to check that $\fX^* \geqdom \bX^*$.  By induction there is a
unique tree $T^*$ with $n$ leaves such that $(\bX(T^*),\fX(T^*)) =
(\bX^*,\fX^*)$. Making the $i$th leaf of $T^*$ an inner vertex with
two leaf children yields the unique tree $T$ with $(\bX(T),\fX(T)) =
(\bX,\fX)$.  \qed

There is a natural correspondence between
$0,1$ strings $\sigma \in \binom{k+\ell}{k}$ and
\term{upright lattice paths}~$P_\sigma$ from $(0,0)$ to $(\ell,k)$,
which takes an entry $1$ from $\sigma$ to a step to the right,
i.e., the addition of $(1,0)$ to the current position, and an
entry $0$ from $\sigma$ to a step upwards,
i.e., the addition of $(0,1)$ to the current position.

This correspondence is the heart of a correspondence between pairs
$(\sigma,\tau) \in \binom{k+\ell}{k}$ with $\tau \geqdom \sigma$, and pairs
$(P_\sigma,P_\tau)$ of non-intersecting lattice paths, where
$P_\sigma$ is from $(0,1)$ to $(k,\ell+1)$ and $P_\tau$ is from $(1,0)$
to $(k+1,\ell)$. This yields a cryptomorphic version of
Theorem~\ref{thm:tree-2strings}.

\begin{theorem}\label{thm:tree-2paths}
  There is a bijection between full binary trees with $k+1$ left leaves
  and $\ell+1$ right leaves and pairs $(P_\bp,P_\fp)$ of
  non-intersecting upright lattice paths, where $P_\bp$ is from
  $(0,1)$ to $(\ell,k+1)$ and $P_\fp$ is from $(1,0)$ to $(\ell+1,k)$.
\end{theorem}

The advantage of working with non-intersecting lattice paths is
that now we can apply the Lemma of
Gessel-Viennot~\cite{gv-bdphlf-85}; see also~\cite{ai-ce-07}.

\begin{theorem}
   The number of full binary trees with $k+1$ left leaves
   and $\ell+1$ right leaves is
   $$
   \det\begin{pmatrix} {k+\ell \choose k} &  {k+\ell \choose k-1} \\[1mm]
                       {k+\ell \choose k+1} &  {k+\ell \choose k}
         \end{pmatrix}
                = \frac{1}{k+\ell+1}{k+\ell+1 \choose k}{k+\ell+1 \choose k+1}
$$
\end{theorem}
This is the \term{Narayana number} $N(k+\ell+1,k+1)$.
From an elementary application of Vandermond's convolution,
$\sum_{k=1}^{n-1} N(n,k) = \frac{1}{n}{2n \choose n-1} = C_n$.
The following proposition summarizes our findings about
Narayana families.

\begin{proposition}
The Narayana number $N(k+\ell+1,k+1)$ counts
\Bitem alternating trees with $k+1$ left vertices and $\ell+1$ right vertices,
\Bitem full binary trees with $k+1$ left leaves and $\ell+1$ right leaves,
\Bitem pairs $(\sigma,\tau)$ of $0,1$ strings in $\binom{k+\ell}{k}$
       with $\tau \geqdom \sigma$,
\Bitem pairs $(P_1,P_2)$ of
  non-intersecting upright lattice paths, where $P_1$ is from
  $(0,1)$ to $(k,\ell+1)$ and $P_2$ is from $(1,0)$ to $(k+1,\ell)$.
\end{proposition}

% *************************************************************
\section{Twin Pairs of Trees and the Baxter Numbers}\label{sec:twins+baxt}
% *************************************************************
\label{sec:paths}

After the Catalan digression we come back to twin pairs of trees.

\begin{definition}
A pair $(A,B)$ of full binary trees
whose fingerprints satisfy
$\fX(A) = \rho(\fX(B))$ is called a
\term{twin-binary pair of trees}.
\end{definition}

\begin{theorem}\label{thm:ta+tb}
There is a bijection between twin-alternating pairs of trees
on $n$ vertices and twin-binary pairs of trees with $n$ leaves.
\end{theorem}

\Proof Let $(A,B)$ be twin-binary trees. Apply the correspondence
from Proposition~\ref{prop:fingerp-bij} to both. This yields trees
$S$ and $T$ such that  $\fprX(S) =
\fX(A)$ and $\fprX(T) = \fX(B)$ and $S^\lambda = A$ and $T^\lambda = B$.
From $\fX(A) = \rho(\fX(B))$ we
conclude $\fprX(S) = \rho(\fprX(T))$ which is the defining
property for twin-alternating trees. \qed

In the proof of Theorem~\ref{thm:2or+ta} we have seen how a
twin-alternating pair of trees can be extended and then glued
together to yield a 2-book embedding of a quadrangulation; see
also Figure~\ref{fig:small-ex}. Doing a similar gluing for a
twin-binary pair of trees, with both trees drawn as in the proof
of Proposition~\ref{prop:fingerp-bij}, yields a particular
rectangulation of the square. Figure~\ref{fig:rectangulation}
shows an example.

%%
%%%%%%%%%%%%%%%%%%%%%%%%%%%%%%%%%%%%%%%%%%%%%%%%%%%
%%
% in einem figure environment mit caption
   \calc_figscale{17}
    \begin{figure}[htb]
    \centerline{\input{\path/rectangulation.pstex_t}}
    \caption{\label{fig:rectangulation}}
    \end{figure}
    VC
{A twin-binary pair of trees and the associated rectangulation.}
%%
%%%%%%%%%%%%%%%%%%%%%%%%%%%%%%%%%%%%%%%%%%%%%%%%%%%
%%

\begin{definition}\label{def:rectang}
Let $X$ be a set of points in the plane and let $R$ be an axis-aligned
rectangle which contains $X$ in its open interior. A
\term{rectangulation of $X$} is a subdivision of $R$ into rectangles
by non-crossing axis-parallel segments, such that every segment
contains a point of~$X$ and every point lies on a segment.
\end{definition}

We are mainly interested in rectangulations of diagonal
sets, i.e., of the sets $X_{n-1} = \{ (i,n-i) : 1\leq i \leq n-2 \}$.
In this case the enclosing rectangle $R$ can be chosen to be the
square spanned by $(0,0)$ and $(n,n)$. Figure~\ref{fig:rectangulation}
shows a rectangulation of $X_{13}$. The following theorem
is immediate from the definitions.

\begin{theorem}\label{thm:tb+rect}
There is a bijection between twin-binary pairs of trees
with $n$ leaves and rectangulations of $X_{n-2}$.
\end{theorem}

\Note
Hartman et al.~\cite{bnz-gig-91} and later independently
de Fraysseix et al.~\cite{fmp-lrspg-95}
prove that it is possible to assign a set of internally disjoint
vertical and horizontal segments to the vertices of any bipartite
graph $G$ such that two segments touch if, and only if, there is an edge between
the corresponding vertices. A proof of this result can be given along
the following line: Extend $G$ by adding edges and vertices to a
quadrangulation $Q$. Augment $Q$ with a 2-orientation and trace the
mappings from 2-orientations via twin-alternating pairs of trees to a
rectangulation of a diagonal point set.  The horizontal and vertical
segments through the points are a touching segment representation for
$Q$. Deleting some and retracting the ends of some other segments
yields a representation for $G$.  A similar observation was made by
Ackerman, Barequet and Pinter~\cite{abp-nrpps-06}.
\EndNote

\ni
Let $(S,T)$ be a twin pair of binary trees with $k+1$ left and
$\ell+1$ right leaves. The bijection from Theorem~\ref{thm:tree-2paths}
maps $Y \in \{S,T\}$ to a pair $(P_\bp(Y),P_\fp(Y))$ of
non-intersecting upright lattice paths, where $P_\bp(Y)$ is from
$(0,1)$ to $(\ell,k+1)$ and $P_\fp(Y)$ is from $(1,0)$ to
$(\ell+1,k)$. Since by definition $\fX(S) = \rho(\fX(T))$
a point reflection of $P_\fp(T)$ at $(0,0)$ followed by a translation
by $(\ell+2,k)$ maps $P_\fp(T)$ to the path $P^*_\fp(T)$ defined as $P^*_\fp(T)=P_\fp(S)$.
The same geometric transformation maps $P_\bp(T)$ to $P^*_\bp(T)$
from $(2,-1)$ to $(\ell+2,k-1)$ and of course $(P^*_\bp(T),P^*_\fp(T))$ is again a pair of
non-intersecting upright lattice paths. Actually,
$(P_\bp(S),P_\fp(S)  ,P^*_\bp(T)) =
 (P_\bp(S),P^*_\fp(T),P^*_\bp(T))$ is a triple of
non-intersecting upright lattice paths. Since the first two of these
paths uniquely determine $S$ and the last two uniquely determine $T$
we obtain, via a translation
of the three paths one unit up, the following theorem.

%%%%%%%%%%%%%%%%%%%%%%%%%%%%%%%%%%%%%%%%%%%%%%%%%%%
%%
% in einem figure environment mit caption
   \calc_figscale{42}
    \begin{figure}[htb]
    \centerline{\input{\path/three-paths.pstex_t}}
    \caption{\label{fig:three-paths}}
    \end{figure}
    VC
{A twin-binary pair of trees and its triple
of non-intersecting lattice paths.}
%%
%%%%%%%%%%%%%%%%%%%%%%%%%%%%%%%%%%%%%%%%%%%%%%%%%%%

\begin{theorem}\label{thm:twins-3paths}
  There is a bijection between twin pairs of full binary trees with $k+1$ left leaves
  and $\ell+1$ right leaves and triples $(P_1,P_2,P_3)$ of
  non-intersecting upright lattice paths, where $P_1$ is from
  $(0,2)$ to $(k,\ell+2)$, $P_2$ is from $(1,1)$ to $(k+1,\ell+1)$,
  and $P_3$ is from $(2,0)$ to $(k+2,\ell)$.
\end{theorem}

Again we can apply the Lemma of Gessel-Viennot.

\begin{theorem}
   The number of twin pairs of full binary trees with $k+1$ left leaves
  and $\ell+1$ right leaves is
   $$
    \det\begin{pmatrix}
        {k+\ell \choose k}   & {k+\ell \choose k-1} & {k+\ell \choose k-2}\\[1mm]
        {k+\ell \choose k+1} & {k+\ell \choose k}   & {k+\ell \choose k-1}\\[1mm]
        {k+\ell \choose k+2} & {k+\ell \choose k+1} & {k+\ell \choose k}
         \end{pmatrix}
= 2\;\frac{(k+\ell)!\; (k+\ell+1)!\; (k+\ell+2)!}%%
{k!\;(k+1)!\; (k+2)! \; \ell! \; (\ell+1)! \; (\ell+2)!} = \Theta_{k,\ell}
$$
\end{theorem}

The number $\Theta_{k,\ell}$ has some quite nice expressions in terms
of binomial coefficients, e.g.,
$\Theta_{k,\ell} = \frac{2}{(k+1)^2\, (k+2)}
       {k+\ell \choose k}{k+\ell+1 \choose k}{k+\ell+2 \choose k}$
or
$\Theta_{k,\ell} = \frac{2}{(n+1)\, (n+2)^2}
       {k+\ell+2 \choose k}{k+\ell+2 \choose k+1}{k+\ell+2 \choose k+2}$.
The total number of twin binary trees with $n+2$ leaves is given by the
\term{Baxter number}
$$
B_{n+1} = \sum_{k=0}^{n} \Theta_{k,n-k},
$$
whose initial values are $1, 2, 6, 22, 92, 422, 2074, 10754$.
The next proposition collects
families that are, due to our bijections, enumerated by $\Theta$-numbers.

\begin{proposition}\label{prop:theta-count}
The number $\Theta_{k,\ell}$ counts
\Bitem triples $(P_1,P_2,P_3)$ of
  non-intersecting upright lattice paths, where $P_1$ is from
  $(0,2)$ to $(\ell,k+2)$ and $P_2$ is from $(1,1)$ to $(\ell+1,k+1)$ and
  $P_3$ is from $(2,0)$ to $(\ell+2,k)$ .
\Bitem twin pairs of binary trees with $k+1$ left leaves and $\ell+1$ right leaves,
\Bitem rectangulations of $X_{k+\ell}$ with $k$ horizontal and $\ell$ vertical segments,
\Bitem twin pairs of alternating trees with $k+1$ left vertices and $\ell+1$ right vertices,
\Bitem separating decompositions of quadrangulations
       with $k+2$ white and $\ell+2$ black vertices,
\Bitem 2-orientations of quadrangulations with $k+2$ white and $\ell+2$ black vertices.
\end{proposition}

\Note
The concept of twin-binary pairs of trees is due to Dulucq and
Guibert~\cite{dg-swstbp-96}. They also give a bijection between
twin-binary pairs of trees and triples of non-intersecting
lattice paths. The bijection also uses the fingerprint as the
middle path, the other two are defined differently.
In~\cite{dg-bp-98} they extend their work to include some
more refined counts. A very good entrance point for more
information about Baxter numbers is The On-Line Encyclopedia
of Integer Sequences~\cite[A001181]{oleis}.

\FreeItem{}
Fusy, Schaeffer and Poulalhon~\cite{fsp-bcbo-06} gave a direct bijection from
separating decompositions to triples of non-intersecting paths in a grid.
Their main application is the counting of bipolar orientations of rooted 2-connected
maps. These results are included in Section~\ref{ssec:bipol}.

\FreeItem{}
Ackerman, Barequet and Pinter~\cite{abp-nrpps-06} also have the result
that the number of rectangulations of~$X_{n}$ is the Baxter number
$B_{n+1}$.  Their proof is via a recurrence formula obtained by Chung
et al.~\cite{cghk-nbp-78}.  They also show that for a point set
$X_{\pi}=\{(i,\pi(i)) : 1\leq i \leq n\}$ to have exactly $B_{n+1}$
rectangulations it is sufficient that $\pi$ is a Schr\"oder
permutation, i.e., a permutation avoiding the patterns $3-1-4-2$ and
$2-4-1-3$.  They conjecture that whenever $\pi$ is a permutation that
is not Schr\"oder, the number of rectangulations of $X_{\pi}$ is
strictly larger than the Baxter number.

\FreeItem{}
In contrast to the nice formulas for the number of 2-orientations of
quadrangulations on $n$ vertices, very little is known about the
number of 2-orientations of a fixed quadrangulation $Q$.
In~\cite{fz-nao-07} it is shown that the maximal number of 2-orientations a
quadrangulation on $n$ vertices can have is asymptotically between
$1,47^n$ and $1,91^n$. To our knowledge, 
the computational complexity of the counting problem is open.
\EndNote

%%%%%%%%%%%%%%%%%%%%%%%%%%%%%%%%%%%%%%%%%%%%%%%%%%%%%%%%%%%%%%%%%%%%%%%%%

% *************************************************************
\section{More Baxter Families}\label{sec:more_baxt}
% *************************************************************

In this section we deal with Baxter permutations and bipolar orientations.
In both families we identify objects counted by $\Theta$-numbers and
less refined families counted by Baxter numbers.

% *************************************************************
\subsection{Baxter Permutations}\label{ssec:baxt_perm}
% *************************************************************

\begin{definition}
The \term{max-tree} $\Max(\pi)$ of a permutation $\pi$ is recursively
defined as the binary tree with root labeled $z$, left
subtree $\Max(\pi_{\text{left}})$ and right subtree
$\Max(\pi_{\text{right}})$ where $z$ is the maximum entry of $\pi$
and in one-line notation $\pi = \pi_{\text{left}},z,\pi_{\text{right}}$.
The recursion ends with unlabeled leaf-nodes corresponding
to  $\Max(\emptyset)$.
\end{definition}

The max-tree of a permutation is a  full binary tree.
The $i$th leaf $v_i$ of $\Max(\pi)$ from the left
corresponds to the adjacent pair $(\pi_{i-1},\pi_i)$ in the
permutation $\pi$. Leaf $v_i$ is a left leaf if, and only if, $(\pi_{i-1},\pi_i)$ is a
descent, i.e., if $\pi_{i-1} > \pi_i$.

The \term{min-tree} $\Min(\pi)$ of a permutation $\pi$ is defined dually,
i.e., as the binary
tree with root labeled $a$, left subtree $\Min(\pi_{\text{left}})$ and
right subtree $\Min(\pi_{\text{right}})$ where $a$ is the minimum entry
of $\pi = \pi_{\text{left}},a,\pi_{\text{right}}$. Also let $\Min(\emptyset)$
be a leaf-node. The $i$th leaf $y_i$ of $\Min(\pi)$ from the left
is a left leaf if, and only if, $(\pi_{i-1},\pi_i)$ is a rise,  i.e., if $\pi_{i-1} < \pi_i$.

With these definitions and observations, see also Figure~\ref{fig:min-max-tree}, we obtain:

\begin{proposition}\label{prop:min-max-trees}
For a permutation $\pi$ of $[n-1]$ the pair
$(\Max(\pi),\Min(\rho(\pi)))$
is a twin binary pair of trees.
\end{proposition}
%%
%%%%%%%%%%%%%%%%%%%%%%%%%%%%%%%%%%%%%%%%%%%%%%%%%%%
%%
% in einem figure environment mit caption
   \calc_figscale{22}
    \begin{figure}[htb]
    \centerline{\input{\path/min-max-tree.pstex_t}}
    \caption{\label{fig:min-max-tree}}
    \end{figure}
    VC
{The trees $\Max(\pi)$ and $\Min(\rho(\pi))$ associated with $\pi = 1,7,4,6,3,2,5$.}
%%
%%%%%%%%%%%%%%%%%%%%%%%%%%%%%%%%%%%%%%%%%%%%%%%%%%%
%%

This mapping from permutations to twin binary pairs of trees is
not a bijection. Indeed, the permutation $\pi' = 1,7,5,6,3,2,4$
also maps to the pair of trees shown in
Figure~\ref{fig:min-max-tree}.

\begin{definition}
A \term{Baxter permutation} is a permutation which avoids the
pattern $2 \tsl 4 1 \tsl 3$ and $3 \tsl 1 4 \tsl 2$. That is,
$\pi$ is Baxter if there are no indices $i < j , j+1 < k$ with
$\pi_{j+1} < \pi_i < \pi_k < \pi_j$ nor with $\pi_{j+1} > \pi_i >
\pi_k > \pi_j$.
\end{definition}

\begin{theorem}\label{thm:tb+Baxt}
There is a bijection between twin binary trees with $n$ leaves and
Baxter permutations of $[n-1]$.
\end{theorem}

\Proof From Proposition~\ref{prop:min-max-trees} we know that
$(\Max(\pi),\Min(\rho(\pi)))$ is a twin binary pair of trees; this
remains true if we restrict $\pi$ to be Baxter.

For the converse let $(S,T)$ be twin binary trees with $n$
leaves. With this pair associate a rectangulation of $X_{n-2}$.
Tilt the rectangulation to get the diagonal points onto the
$x$-axis. Each rectangle of the rectangulation contains a
highest corner which corresponds to an inner vertex of $S$
and a lowest corner corresponding to an inner vertex of $T$.
We will refer to these corners as the \term{north-corner}
and the \term{south-corner} respectively.

The idea is to associate a number with every rectangle; the
permutation $\pi$ corresponding to $(S,T)$ is then read of from
the order of intersection of the rectangles with the $x$-axis.
Writing the numbers of rectangles to their north- and
south-corners makes $(S,T)= (\Max(\pi),\Min(\rho(\pi)))$.

The algorithm associates the numbers with rectangles in decreasing
order. Number $n-1$ is associated with the rectangle with highest
north-corner. After having associated $k$ with some rectangle
$R_k$, the union of unlabeled rectangles can be seen as a series
of pyramids over the $x$-axis; see Figure~\ref{fig:pyramids}.
%%
%%%%%%%%%%%%%%%%%%%%%%%%%%%%%%%%%%%%%%%%%%%%%%%%%%%
%%
% in einem figure environment mit caption
   \calc_figscale{15}
    \begin{figure}[htb]
    \centerline{\input{\path/pyramids.pstex_t}}
    \caption{\label{fig:pyramids}}
    \end{figure}
    VC
{A series of pyramids after labeling rectangle $k$.}%%
%%
%%%%%%%%%%%%%%%%%%%%%%%%%%%%%%%%%%%%%%%%%%%%%%%%%%%

The label $k-1$ has to correspond to one of the rectangles which
have their north-corner on the tip of one of the pyramids (this is
because $S$ will be the max-tree of $\pi$). The algorithm will
choose the next pyramid to the left or the next pyramid to the
right of the interval on the $x$-axis which belongs to $R_k$. The
decision which of the two is taken depends on the south-corner of
$R_k$. If the south-corner is a
\lower3pt\hbox{\includegraphics{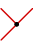}}, i.e., a left child in
tree $T$, then the pyramid to the left is chosen; otherwise, if
the south-corner is a
\lower3pt\hbox{\includegraphics{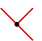}}\kern1.5pt, then the
pyramid to the right is chosen.

%%
%%%%%%%%%%%%%%%%%%%%%%%%%%%%%%%%%%%%%%%%%%%%%%%%%%%
%%
% in einem figure environment mit caption
   \calc_figscale{27}
    \begin{figure}[htb]
    \centerline{\input{\path/square-labeling.pstex_t}}
    \caption{\label{fig:square-labeling}}
    \end{figure}
    VC
{Generating a Baxter permutation from a rectangulation.
 The final state shows the permutation with its max- and min-trees.}%%
%%
%%%%%%%%%%%%%%%%%%%%%%%%%%%%%%%%%%%%%%%%%%%%%%%%%%%
%%
An example for the execution of this algorithm is given in
Figure~\ref{fig:square-labeling}.
The proof that this is a bijection is done with three claims.

\Claim A. The permutation constructed by the algorithm is Baxter.
\medskip

\ni{\it Proof of the claim.}
Let $a<b<c<d$, we want to show that the algorithm will not produce
the pattern $b \tsl d a \tsl c$. Think of the status after labeling
rectangle $R(c)$. To have a chance of producing the pattern $d$ is
left of $c$ and the slot immediately to the right of $d$ is not yet
used, i.e., it belongs to a pyramid $P$. From the labeling rule of the algorithm
it follows that the rectangle covering the leftmost slot of $P$
has to be labeled before any rectangle left of $P$ can be labeled.
This shows that the pattern is impossible. The case of the other
pattern is symmetric.\qedclaim

\Claim B. The border between labeled and unlabeled rectangles at any
stage of the algorithm is a zig-zag (the definition of zig-zag should
be evident from Figure~\ref{fig:pyramids}).\medskip

%%%%%%%%%%%%%%%%%%%%%%%%%%%%%%%%%%%%%%%%%%%%%%%%%%%
%%
% in einem figure environment mit caption
   \calc_figscale{13}
    \begin{figure}[htb]
    \centerline{\input{\path/no-zigzag.pstex_t}}
    \caption{\label{fig:no-zigzag}}
    \end{figure}
    VC
{A violation of the zig-zag property.}%%
%%
%%%%%%%%%%%%%%%%%%%%%%%%%%%%%%%%%%%%%%%%%%%%%%%%%%%

\ni{\it Proof of the claim.} Suppose not. Then there is a first
rectangle $R$ whose labeling violates the property,
Figure~\ref{fig:no-zigzag} shows the situation up to a reflection
interchanging left and right. Let~$y$ be the label of $R$ and let
$z$ be the label of the rectangle whose south-corner is the
deepest point of the valley whose shape was destroyed by $R(y)$.
The rule of the algorithm implies that the rectangle labeled after
$z$ was left of $R(z)$. This rectangle $R(z-1)$ has to intersect
the $x$-axis somewhere left of pyramid $P$. From the labeling rule
of the algorithm it follows that the rectangle covering the
rightmost slot of $P$ has to be labeled before any rectangle right
of~$P$ can be labeled. Rectangle $R$ however is right of $P$, a
contradiction.\qedclaim

Claim B implies that in the lower tree the labels of south-corners of
rectangles are decreasing along every path from a leaf to the root.
This is a property characterizing $\Min$-trees, hence, the labels of
the rectangles yield the $\Min(\rho(\pi))$. The fact that the
labeling of the north-corners of rectangles yields the $\Max(\pi)$ is
more immediate from the algorithm.

\Claim C. If $\sigma$ is a permutation with
$(\Max(\sigma),\Min(\rho(\sigma))) = (S,T)$ and $\sigma$ is not the result
of applying the algorithm to the rectangulation corresponding to
$(S,T)$, then $\sigma$ is not Baxter.
\medskip

\ni{\it Proof of the claim.} Compare the one-line notation of
$\sigma$ and $\pi$, where $\pi$ is the result of applying the
algorithm to the rectangulation corresponding to $(S,T)$. Consider
the largest value $k$ which is not at the same position in
$\sigma$ and $\pi$. Clearly $k\neq n-1$, because $\sigma$ and
$\pi$ have the same max-tree. Hence, if we think of the algorithm
producing $\pi$ in the state when $k+1$ is placed, then the
placement of $k$ in $\sigma$ fails to obey the rules of the
algorithm. Either $k$ is placed as to have its north-corner in a
pyramid on the wrong side, or the side is respected but the
pyramid containing the north-corner of $k$ is not next to $R(k+1)$
on this side. We indicate how to find a forbidden pattern in each
of the two cases.

{\bf Wrong side.}
Suppose the south-corner $p$ of $R(k+1)$ is a
\lower3pt\hbox{\includegraphics{y.eps}} and $k$ is placed to the right
of $k+1$. Let $a$ be the element in $\sigma$ whose rectangle has its
east-corner at $p$ and let $q$ be the first node of type
\lower3pt\hbox{\includegraphics{co-y.eps}} on the path from $p$ to the
root of the min-tree, i.e, of $T$. Let $b$ be the element in $\sigma$
whose rectangle has its west-corner at $q$. From the
min-tree property we  infer that $b<a<k$. Since $k+1$ and $b$ are
neighbors in $\sigma$, the elements $a,k+1,b,k$ form a forbidden $3\tsl 1 4
\tsl 2$ pattern. 

The other case where the south-corner $p$ of $R(k+1)$ is a
\lower3pt\hbox{\includegraphics{co-y.eps}} is symmetric. In this case
there is a forbidden pattern $2 \tsl 4 1 \tsl 3$.

{\bf Wrong pyramid.}
Suppose the south-corner $p$ of $R(k+1)$ is a
\lower3pt\hbox{\includegraphics{y.eps}} and $k$ is placed to a pyramid
left of $k+1$ but not to the first one. Let $r>k+1$ be some element
separating the first from the second pyramid.
Let $a$ be the element in $\sigma$ whose rectangle has its
east-corner at $p$ and note that $a<k$. Between $r$ and $a$ there is
an adjacent pair $r',a'$ with $a' < k$ and $k+1 < r'$. Hence,
$k,r',a',k+1$ is a forbidden  $2 \tsl 4 1 \tsl 3$. Again, the second
case is symmetric. This completes the proof of the claim
\qedclaim

\ni
Claim A says that every twin binary pair $(S,T)$ of trees
is mapped by the algorithm to a Baxter permutation $\pi$.
As a consequence of Claim B we noted that $(S,T) = (\Max(\pi),\Min(\rho(\pi))$.
Claim C is the injectivity, hence, the mapping is a
bijection.
\qed

From Proposition~\ref{prop:theta-count} and the observation about the correspondence
of left and right leaves in the max-tree of a permutation with descents and rises
we obtain:

\begin{proposition}\label{prop:Mallows}
The number $\Theta_{k,\ell}$ counts
\Bitem twin pairs of binary trees with $k+1$ left leaves and $\ell+1$ right leaves,
\Bitem Baxter permutations of $k+\ell+1$ with $k$ descents and $\ell$ rises.
\end{proposition}

\Note
Baxter numbers first appeared in the context of counting Baxter permutations.
Chung, Graham, Hoggatt and Kleiman~\cite{cghk-nbp-78} found some interesting recurrences and gave a proof based on generating functions. Mallows~\cite{m-bpra-79}
found the refined count of Baxter permutations by rises (Proposition~\ref{prop:Mallows}).
The bijection of Theorem~\ref{thm:tb+Baxt} is essentially due to
Dulucq and Guibert~\cite{dg-swstbp-96,dg-bp-98}. Their description and proof,
however, does not use geometry.
They also prove Proposition~\ref{prop:Mallows} and
some even more refined counts, e.g., the number of Baxter permutations of $[n]$ with
$\ell$ rises and $s$ left-to-right maxima and $t$ right-to-left maxima.
\EndNote

\ni
A permutation $(a_1,a_2,\ldots,a_n)$ is \term{alternating} if
$a_1 < a_2 > a_3 < a_4 > \ldots $, i.e., each consecutive
pair $a_{2i-1}, a_{2i}$ is a rise and each pair $a_{2i}, a_{2i+1}$
a descent.  Alternating permutations are characterized
by the property that the reduced fingerprints of their $\Min$- and $\Max$-trees
are alternating, i.e., of the form $\ldots0,1,0,1,0,1,\ldots$
and in addition, to ensure that the first pair is a rise,
the first entry of the reduced fingerprints of the $\Max$-tree is a 0.
Due to this characterization we obtain the following specialization of
Theorem~\ref{thm:tb+Baxt}:

\begin{lemma}
Twin pairs of binary trees with
an alternating reduced fingerprint starting in 0 and alternating Baxter permutations
are in bijection.
\end{lemma}

Let $T$ be a binary tree with $n$ leaves and with an alternating reduced fingerprint starting
with a 0. The leaves of $T$ come in pairs from left to right so that
the leaves from each pair are attached to the same interior node.
Pruning the leaves we obtain a tree $T'$ with $n-\lfloor\frac{n}{2}\rfloor$ leaves.
From $T'$ we come back to $T$ by attaching a new pair of leaves to each of the first
$\lfloor\frac{n}{2}\rfloor$ leaves of~$T'$. Using this kind of bijection we obtain
two bijections (see Figure~\ref{fig:alt-baxt}):
\Bitem
a bijection between alternating Baxter permutations of $[2k-1]$ and pairs of
binary trees with $k$ leaves, and
\Bitem
a bijection between alternating Baxter permutations of $[2k]$ and pairs of
binary trees with $k$ and $k+1$ leaves.

\begin{theorem}\label{thm:alt-baxt}
The number of alternating Baxter permutations on $[n-1]$ is
$C_{k-1}C_{k}$ if $n=2k$ and $C_{k-1}C_{k-1}$ if $n=2k-1$.
\end{theorem}

%%%%%%%%%%%%%%%%%%%%%%%%%%%%%%%%%%%%%%%%%%%%%%%%%%%
%%
% in einem figure environment mit caption
   \calc_figscale{33}
    \begin{figure}[htb]
    \centerline{\input{\path/alt-baxt.pstex_t}}
    \caption{\label{fig:alt-baxt}}
    \end{figure}
    VC
{Alternating Baxter permutations and pairs of trees.}%%
%%
%%%%%%%%%%%%%%%%%%%%%%%%%%%%%%%%%%%%%%%%%%%%%%%%%%%

\vskip-5mm
\Note
Theorem~\ref{thm:alt-baxt} was obtained by Cori et al.~\cite{cdv-spsbp-86}.
It was reproved by Dulucq and Guibert~\cite{dg-swstbp-96} as a specialization
of their bijection between Baxter permutations and twin pairs of binary trees.
In~\cite{gl-dabpac-00} it is shown that alternating Baxter permutations with the property
that their inverse is again alternating Baxter are counted by the Catalan numbers.
\EndNote

% *************************************************************
\subsection{Plane Bipolar Orientations}\label{ssec:bipol}
% *************************************************************
A graph $G$ is said to be \emph{rooted} if one of its edges is distinguished and oriented.
The origin and the end of the root-edge are denoted  $s$ and $t$. If $G$ is a \emph{plane graph}, the root-edge is always assumed to be incident to the outer face,
with the outer face on its left.

\begin{definition}
A \term{bipolar orientation} of a rooted graph $G$  is an acyclic
orientation of $G$ such that the unique source (i.e., vertex with only outgoing edges) is $s$ and the unique sink (i.e., vertex with only ingoing edges) is $t$.
A \emph{plane bipolar orientation} is a bipolar orientation on a rooted plane graph
(multiple edges are allowed).
\end{definition}

%\ni
%We consider bipolar orientations of plane graphs with specified poles $s$ and $t$
%on the outer face. In most cases we additionally assume that $s$ and $t$ are
%neighbors on the outer face; in that case we say that the graph is
%\term{rooted} at the edge $(s,t)$.
It is well-known that the rooted graphs admitting a bipolar
orientation are exactly 2-connected graphs,
i.e., graphs with no separating vertex.

\Note
 Bipolar orientations have proved to be insightful in solving
many algorithmic problems such as planar graph
embedding~\cite{Le66,Chiba} and geometric representations of
graphs in various flavors (visibility~\cite{TaTo}, floor
planning~\cite{Ro,KaHe}, straight-line drawing~\cite{TaTo2,Fu6}).
They also constitute a beautiful combinatorial structure; the
thesis of Ossona de Mendez is devoted to studying their numerous
properties and applications~\cite{OssThes}; see also~\cite{DeOss}
for a detailed survey.
\EndNote

%Let a \term{rooted map} be a rooted 2-connected plane multigraph $M$.
%Given a rooted map, let $\theta(M)$ be the number of
%bipolar orientations of $M$.
%The purpose of this subsection is to show an explicit formula for the
%sum of $\theta(M)$ where $M$ is running through all maps with the same size
%parameters.  Precisely, let $\cM_{ij}$ be the set of rooted
%maps with $(i+1)$ vertices and $(j+1)$ faces
%(including the outer one); then the
%quantity $\sum_{M\in\cM_{ij}}\theta(M)$ turns out to be $\Theta_{ij}$
%(Proposition~\ref{prop:theta-bip}).

%It is known that for a 2-connected
%graph $G$, $\beta = 1$ if and only if  $G$ is series-parallel
%\cite{bry-ox-92}. Hence series-parallel graphs (or maps) are the
%only 2-connected ones having a unique bipolar orientation up to
%reversion. 

%\begin{definition}\label{angular}
Let $G$ be a rooted plane graph; the \term{angular map} of $G$
is the graph $Q$ with vertex set consisting of vertices and faces of $G$,
and edges corresponding to incidences between a vertex and a face.
The special vertices $s,t$ of $Q$ are the extremities (origin $s$ and end $t$)
of the root-edge of $G$.

%\end{definition}

The angular map $Q$ of $G$ inherits a plane embedding from $G$.
The unique bipartition of~$Q$ has the vertices of $G$ in one color
class and the faces of $G$ in the other. We assume that vertices
of $G$ are black and faces of $G$ are white. All the faces of $Q$
are quadrangles, which correspond to the edges of $G$.
Moreover, since  $G$ is 2-connected, $Q$ has no double edges, so
$Q$ is a quadrangulation. 
%It is well known that the angular map is a
%bijection between rooted 2-connected plane graphs and quadrangulations.

%%%%%%%%%%%%%%%%%%%%%%%%%%%%%%%%%%%%%%%%%%%%%%%%%%%%%%%%%%%%%%%%%%%%%%
% in einem figure environment mit caption
   \calc_figscale{17}
    \begin{figure}[htb]
    \centerline{\input{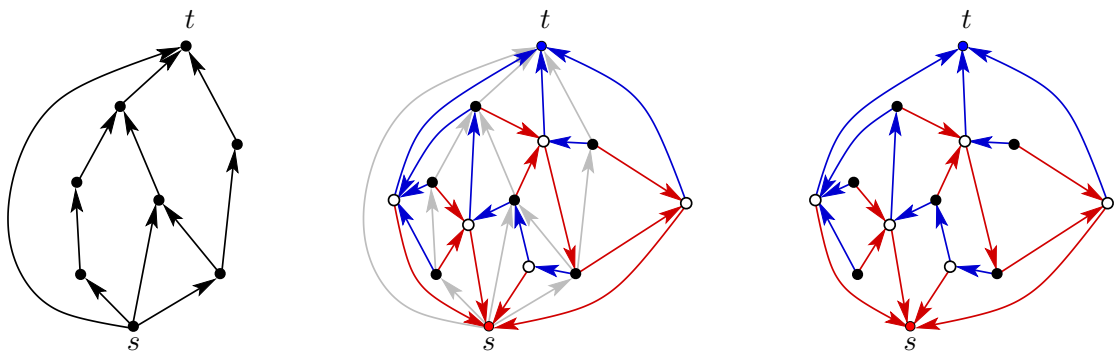}}
    \caption{From a rooted map endowed with a bipolar orientation to
      a separating decomposition on the angular map.\label{fig:bicolor}}
    \end{figure}
    
%%%%%%%%%%%%%%%%%%%%%%%%%%%%%%%%%%%%%%%%%%%%%%%%%%%%%%%%%%%%%%%%%%%%%%

If $G$ is endowed with a bipolar orientation $B$, the angular map
can be enriched in order to transfer the orientation
onto~$Q$. Actually, we will define a bijection between bipolar orientations of $G$ and
separating decompositions of $Q$.
The construction, based on two facts about bipolar orientations of rooted plane graphs,
is illustrated in Figure~\ref{fig:bicolor}.

\Fact V. Every vertex $v\neq s,t$ has exactly two adjacent faces
         (angles) where the orientation of the edges differ.

\Fact F. Every face $f$ has exactly two vertices (angles)
         where the orientation of the edges coincide.

\begin{figure}[htb]
    \centerline{\input{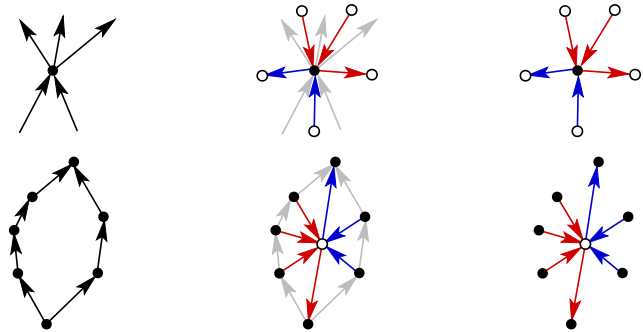}}
    \caption{The transformation for a vertex and a face of a rooted map.\label{fig:bic-detail}}
    \end{figure}

\smallskip \ni
Facts V and F specify two distinguished edges in the angular map for every
non-special vertex and every face. Since every edge of $Q$ is distinguished
either for a vertex or for a face this yields a 2-orientation.
Figure~\ref{fig:bic-detail} indicates how to color this 2-orientation
to get a separating decomposition on $Q$.
From a separating decomposition on $Q$, the unique bipolar orientation
on $G$ inducing $Q$ is easily recovered. (If $Q$ has white vertices of degree 2,
then there are multi-edges on $G$). To summarize:

%%%%%%%%%%%%%%%%%%%%%%%%%%%%%%%%%%%%%%%%%%%%%%%%%%%%%%%%%%%%%%%%%%%%%%
% in einem figure environment mit caption
   \calc_figscale{16}
    
%%%%%%%%%%%%%%%%%%%%%%%%%%%%%%%%%%%%%%%%%%%%%%%%%%%%%%%%%%%%%%%%%%%%%%

\begin{proposition}\label{thm:bi}
Plane bipolar orientations  with $\ell+2$ vertices and $k+2$ faces are
in bijection with separating decompositions of
quadrangulations with $\ell+2$ black vertices and $k+2$ white vertices.
Consequently, the number $\Theta_{k,\ell}$ counts:
\Bitem separating decompositions of quadrangulations
       with $k+2$ white and $\ell+2$ black vertices,
\Bitem Bipolar orientations of rooted plane graphs with
       $k+2$ faces and $\ell+2$ vertices.
\end{proposition}

%EF0303: I have merged the theorem and the proposition into a single proposition
%\begin{proposition}\label{prop:theta-bip}
%The number $\Theta_{k,\ell}$ counts
%\Bitem separating decompositions of quadrangulations
%       with $k+2$ white and $\ell+2$ black vertices,
%\Bitem Bipolar orientations of rooted plane graphs with
%       $k+2$ faces and $\ell+2$ vertices.
%\end{proposition}

In Section~\ref{sec:schnyder} we will use this and some previous bijections to give an
independent proof for a beautiful formula of Bonichon~\cite{b-brmpgpdp-05}
for the number of Schnyder woods on triangulations with $n$ vertices.

\Note
The two facts {\bf V} and {\bf F} have been rediscovered frequently, they can be found, e.g.,
in~\cite{DeOss,Ro,TaTo}. Actually, plane bipolar orientations can be defined via
properties {\bf V} and {\bf F}. The bijection of Proposition~\ref{thm:bi}
is a direct extension of ~\cite[Theo 5.3]{DeOss}. In 2001, 
R.~Baxter~\cite[Eq 5.3]{baxter} guessed that plane bipolar orientations are
counted by the $\Theta$-numbers. His verification
is based on algebraic manipulations on generating functions of plane graphs
weighted by their Tutte polynomials. A simpler proof for this fact was
obtained by Fusy et al.~\cite{fsp-bcbo-06} via a direct bijection 
from separating decompositions to triples of non-intersecting lattice paths.
 Their  bijection presents significant differences from the one presented
 in this article, even if the classes in correspondence are the same.
 The main difference is that they do not treat the blue tree and the red tree
 of a separating decomposition in 
 a symmetric way, as we do here, and their correspondence is less geometric. 
(In their bijection, the blue tree is encoded as a refined Dyck word,  
while  the red edges are encoded  as the sequence of degrees in red of the white vertices.)

\FreeItem{}
 As shown in~\cite{DeOss}, the number of bipolar orientations
of a fixed rooted graph $G$ is equal to twice the
coefficient $[x]T_G(x,y)$ in the Tutte polynomial of $G$. This
coefficient is called Crapo's $\beta$ invariant and it is \#P-hard
to compute \cite{annan-95}. 
To our knowledge, 
the computational complexity of the counting problem restricted to rooted
\emph{plane} graphs is open. (By the angular map bijection, 
it is clearly equivalent to computing
the number of 2-orientations of a fixed quadrangulation.) 
 \EndNote

%%%%%%%%%%%%%%%%%%%%%%%%%%%%%%%%%%%%%%%%%%%%%%%%%%%%%%%%%%%%%%%%%%%%%%
\subsection{Digression: Duality, Completion Graph, and Hamiltonicity}\label{ssec:hamiltonicity}
%%%%%%%%%%%%%%%%%%%%%%%%%%%%%%%%%%%%%%%%%%%%%%%%%%%%%%%%%%%%%%%%%%%%%%

There exists a well-known \term{duality mapping} for plane graphs. The dual $G^*$
of a plane graph~$G$  has its vertices corresponding to the faces of $G$,
and has its edges corresponding to the adjacencies of the faces of $G$ (two faces are adjacent
if they share an edge). Precisely, each edge $e$ of $G$ gives rise to an edge $e^*$ of $G^*$ that
connects the vertices of $G^*$ corresponding to the faces of $G$ on each side of $e$.
Let us mention here, even if we do not make use of this fact,
that duality can be enriched to take account of bipolar orientations~\cite{DeOss};
if $G$ is endowed with a bipolar orientation, each (oriented) non-root edge $e$ of $G$ gives rise
to an oriented edge $e^*$ of $G^*$ that goes from the face on the left to the face on the right of $e$
 (for the root edge the opposite rule has to be applied).

 %%%%%%%%%%%%%%%%%%%%%%%%%%%%%%%%%%%%%%%%%%%%%%%%%%%%%%%%%%%%%%%%%%%%%%
% in einem figure environment mit caption
   \calc_figscale{19}
    \begin{figure}[htb]
    \centerline{\input{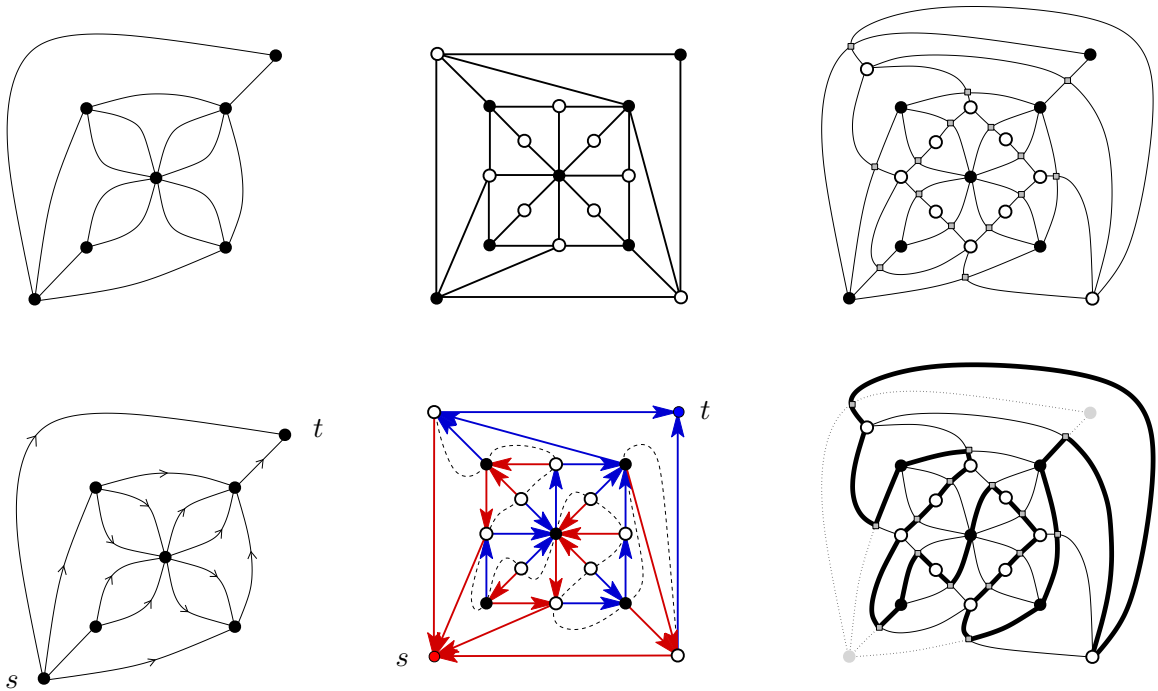}}
    \caption{Top, from left to right: A 2-connected plane graph, its quadrangulation,
and its completion graph. Bottom, from left to right: the same 2-connected plane graph rooted at an
edge and endowed with a bipolar orientation, the quadrangulation endowed with the corresponding
separating decomposition (the equatorial line of the 2-book embedding is drawn ---in dashed line---
by bisecting all bicolored angles),
and the special completion graph endowed with the corresponding Hamiltonian cycle.\label{fig:hamiltonianCyc}}
    \end{figure}
    
%%%%%%%%%%%%%%%%%%%%%%%%%%%%%%%%%%%%%%%%%%%%%%%%%%%%%%%%%%%%%%%%%%%%%%

The \term{completion graph} $\cG$  of $G$ is the plane graph obtained by superimposing
$G$ and its dual $G^*$, see Figure~\ref{fig:hamiltonianCyc}.
Vertices of $\cG$ are of 3 types: \term{primal vertices} are the
vertices of $G$, \term{dual vertices} are the vertices of $G^*$, and \term{edge-vertices}
are the vertices at the intersection of an edge $e\in G$ with its dual edge $e^*$ (hence
edge-vertices have degree 4 in $\cG$). Observe in Figure~\ref{fig:hamiltonianCyc}
that the completion graph of $G$
and the quadrangulation $Q$ of $G$ differ only upon replacing the contour of each face $f$ of $Q$
by a 4-star, the extremities of the 4-star being the 4 vertices incident to $f$ and the center of the 4-star
being the edge-vertex of $\cG$  associated with~$f$ (each edge of $G$ corresponds
both to a face of $Q$ and to an edge-vertex of $\cG$, which yields a correspondence between
faces of $Q$ and edge-vertices of $\cG$).

If $G$ is rooted, the origin $s$ and end $t$ of the root-edge 
are the \term{special vertices}
of~$G$; the \term{special completion graph} of $G$ is the plane graph obtained from 
$\cG$
by removing $s$ and $t$ as well as their incident edges.

\begin{proposition}
The special completion graph of a rooted 2-connected plane graph is Hamiltonian.
\end{proposition}
\Proof
Given $G$ a rooted 2-connected plane graph, we first endow $G$ with a bipolar orientation,
and consider the quadrangulation $Q$ of $G$ endowed with the corresponding separating
decomposition $S$. As proved in  Lemma~\ref{lem:once_faces}), the equatorial
line of $S$
is a simple path that has the two outer nonspecial vertices as extremities and passes once by each inner vertex and
once by the interior of each inner face of $Q$. By the above
discussion on the correspondence between $Q$ and $\cG$, the path
is easily deformed to a simple path on the graph $\cG$ visiting
once all the vertices of $\cG$ except the two special ones, see
Figure~\ref{fig:hamiltonianCyc}. In addition, the path does not
use the 4-star of $\cG$ corresponding to the outer face of $f$.
Completing the path with the two edges of the outer 4-star
incident to white vertices, we obtain a Hamiltonian cycle of the
special completion graph. \qed

% *************************************************************
\section{Symmetries}\label{sec:symm}
% *************************************************************

As we explain in this section, the bijections we have presented have the nice property
that they commute with the
half-turn rotation, which makes it possible to count \emph{symmetric} combinatorial structures.
The first structures we have encountered are 2-orientations. Given a 2-orientation $O$,
exchanging the two special vertices $\{s,t\}$ of $O$ clearly yields another 2-orientation,
which we call the \term{pole-inverted} 2-orientation of $O$ and denote by $\iota(O)$.
A 2-orientation is called \term{pole-symmetric} if $O$ and $\iota(O)$ are isomorphic.

Considering the associated separating decomposition, the blue tree of $O$ is the red tree
of $\iota(O)$ and vice versa. Accordingly, a 2-orientation is pole-symmetric if, and only if, the blue tree
and the red tree are isomorphic as rooted trees, in which case the separating decomposition is
called pole-symmetric as well.  Such a symmetry translates to half-turn rotation symmetries
on the associated embeddings. Indeed, as the two trees composing the separating decomposition
are isomorphic, so are their alternating embeddings and so are the two binary trees that compose
the associated twin pair of full binary trees,
in which case the twin pair is called \emph{symmetric}.
Hence, the 2-book embedding and rectangulation
associated with the separating decomposition are stable under the half-turn rotation that
exchanges the two special vertices. Such 2-book embeddings and rectangulations are called
\emph{pole-symmetric} as well.

Considering  Baxter permutations, the min-tree (resp., max-tree) of a Baxter permutation $\pi$
is the max-tree (resp., min-tree) of the associated Baxter permutation $\overline{\rho(\pi)}$,
i.e., the permutation whose 0-1 matrix is the 0-1 matrix of $\pi$ after half-turn rotation. Baxter
permutations for which $\pi=\overline{\rho(\pi)}$ are said to be \emph{symmetric}.
By definition of the bijective correspondence of Theorem~\ref{thm:tb+Baxt}, a Baxter permutation is
symmetric if, and only if, the associated twin pair of full binary trees is symmetric.

% in einem figure environment mit caption
   \calc_figscale{28}
    \begin{figure}[htb]
    \centerline{\input{\path/symmetriesFig.pstex_t}}
    \caption{\label{fig:symmetriesFig}}
    \end{figure}
    VC
{A pole-symmetric separating decomposition and the corresponding symmetric combinatorial structures: 2-book embedding, twin pair of full binary trees, plane bipolar orientation, Baxter permutation, triple of paths.}

Next we turn to the encoding by a triple of paths. Recall that,
in a twin pair $(T,T')$ of full binary trees, the reduced fingerprints satisfy
 the relation $\hat{\alpha}(T)=\rho(\hat{\alpha}(T'))$. Hence, a symmetric twin pair $(T,T)$
 is characterized by the property that the reduced fingerprint of $T$
satisfies $\hat{\alpha}(T)=\rho(\hat{\alpha}(T))$, i.e., $\hat{\alpha}$ is a palindrome.
Equivalently, if $T$ has $k+1$ left leaves
 and $\ell+1$ right leaves, the upright lattice path $P_2:=P_{\alpha}(T)$, as defined in Section~\ref{sec:paths},
 is stable under the point-reflection $\pi_S$ at $S:=(k/2+1,\ell/2+1)$. The other two paths in the triple
 $(P_1,P_2,P_3)$ of
 non-intersecting lattice paths correspond to two copies of the bodyprint of $T$ read respectively from
 $(0,2)$ to $(k,\ell+2)$ for $P_1$ and from $(\ell+2,k)$ to $(2,0)$ for $P_3$.  Therefore the whole
 triple $(P_1,P_2,P_3)$ is stable under the point reflection $\pi_S$. Such a triple of paths is called
 \emph{symmetric}.

\begin{lemma}
Let  $\thetasym$ be the number of symmetric non-intersecting triples  of upright lattice paths
$(P_1,P_2,P_3)$
going respectively from $(0,2)$, $(1,1)$, $(2,0)$ to $(k,\ell+2)$, $(k+1,\ell+1)$, $(k+2,\ell)$.

(i) If $k$ and $\ell$ are odd, then $\thetasym=0$.

(ii) If $k$ and $\ell$ are even, $k=2\kappa$, $\ell=2\lambda$, then
$$
\thetasym=\sum_{r\geq 1}\frac{2r^3}{(\kappa+\lambda+1)(\kappa+\lambda+2)^2}{\kappa+\lambda+2\choose \kappa+1}{\kappa+\lambda+2\choose \kappa-r+1}{\kappa+\lambda+2\choose \kappa+r+1}.
$$

(iii) If $k$ is odd and $\ell$ is even, $k=2\kappa+1$, $\ell=2\lambda$, then
$$\thetasym=\sum_{r\geq 1}\frac{2r^3\!+\!(\lambda\!-\!r\!+\!1)r(r\!+\!1)(2r\!+\!1)}{(\kappa+\lambda+1)(\kappa+\lambda+2)^2}{\kappa+\lambda+2\choose \kappa+1}{\kappa+\lambda+2\choose \kappa-r+1}{\kappa+\lambda+2\choose \kappa+r+1}.$$

(iv) If $k$ is even and $\ell$ is odd, $k=2\kappa$, $\ell=2\lambda+1$, then
$$\thetasym=\sum_{r\geq 1}\frac{2r^3\!+\!(\kappa\!-\!r\!+\!1)r(r\!+\!1)(2r\!+\!1)}{(\kappa+\lambda+1)(\kappa+\lambda+2)^2}{\kappa+\lambda+2\choose \kappa+1}{\kappa+\lambda+2\choose \kappa-r+1}{\kappa+\lambda+2\choose \kappa+r+1}.$$

\end{lemma}

\Proof By definition, $(P_1, P_2, P_3)$ is stable under the
point-reflection $\pi_S$ at $S:=(\ell/2+1,k/2+1)$.  In particular,
$P_2$ is stable under $\pi_S$, so that $P_2$ has to pass by $S$.  This
can only occur if $S$ is on an axis-coordinate, i.e., $k/2$ or
$\ell/2$ are integers. Therefore $\thetasym=0$ if both $k$ and $\ell$
are odd.

If $k$ and $\ell$ are even, $k=2\kappa$ and $\ell=2\lambda$, the
half-turn symmetry ensures that $P_1,P_2,P_3$ is completely
encoded upon keeping the part $P_1', P_2', P_3'$ of the paths that
lie in the half-plane $\{x+y\leq x_S+y_S\}$, i.e., the half-plane
$\{x+y\leq \kappa+\lambda+2\}$.  The conditions on $(P_1,P_2,P_3)$
translate to the following conditions on the reduced triple:
$(P_1',P_2',P_3')$ is non-intersecting, has same starting points
as $(P_1,P_2,P_3)$, the endpoint of $P_2'$ is $S$, and the
endpoints of $P_1'$ and $P_3'$ are equidistant from $S$, i.e.,
there exists an integer $r\geq 1$ such that $P_1'$ ends at
$(\kappa+1-r,\lambda+1+r)$ and $P_3'$ ends at
$(\kappa+1+r,\lambda+1-r)$.  Hence, up to fixing $r\geq 1$,
$(P_1',P_2',P_3')$ form a non-intersecting triple with explicit
fixed endpoints, so that the number of such triples can be
expressed using Gessel-Viennot determinant formula. The expression
for $\thetasym$ follows.

If $k$ is odd and $\ell$ is even, $k=2\kappa+1$ and $\ell=2\lambda$,
the triple $(P_1,P_2,P_3)$ is again completely encoded by keeping the
part $(P_1',P_2',P_3')$ of the paths that lie in $\{x+y\leq
x_S+y_S\}$, i.e., the half-plane $\{x+y\leq \kappa+\lambda+5/2\}$. The
difference with the case where $k$ and $\ell$ are even is that
$P_1',P_2',P_3'$ are not standard lattice paths, as they end with a
step of length $1/2$.  Similarly as before, the conditions on
$(P_1,P_2,P_3)$ are equivalent to the properties that
$(P_1',P_2',P_3')$ are non-intersecting, have the same starting points
as $(P_1,P_2,P_3)$, $P_2'$ ends at $S$, and $P_1',P_3'$ end at points
that are equidistant from $S$ on the line $\{x+y=x_S+y_S\}$ and have
one integer coordinate, i.e., there exists an integer $m\geq 2$ such
that $P_1'$ ends at $(x_S-m/2,y_S+m/2)$ and $P_3'$ ends at
$(x_S+m/2,y_S-m/2)$.  Notice that, upon discarding the last step, the
system $(P_1',P_2',P_3')$ is equivalent to a triple of
non-intersecting upright lattice paths
$(\overline{P_1'},\overline{P_2'},\overline{P_3'})$ with starting
points $(0,2)$, $(1,1)$, $(2,0)$, and endpoints that are either of the
form $(\kappa+1-r,\lambda+1+r)$, $(\kappa+1,\lambda+1)$,
$(\kappa+1+r,\lambda+1-r)$ if $m$ is even, $m=2r$, or are of the form
$(\kappa+1-r,\lambda+1+r)$, $(\kappa+1,\lambda+1)$,
$(\kappa+2+r,\lambda-r)$ if $m$ is odd, $m=2r+1$.  In each case, the
number of triples has an explicit form from the formula of
Gessel-Viennot.  The expression of $\thetasym$ follows. Finally,
notice that the set of symmetric non-intersecting triples is stable
under swapping $x$-coordinates and $y$-coordinates, yielding the
relation $\thetasym=\Theta_{\ell,k}^{\circlearrowleft}$.  Thus the
formula for $\thetasym$ when $k$ is even and $\ell$ is odd simply
follows from the formula obtained when $k$ is odd and $\ell$ is even.
\qed

Considering bipolar orientations, the effect of the half-turn symmetry
of a separating decomposition on the associated plane bipolar
orientation is clearly that the orientation is unchanged when the
poles are exchanged, the directions of all edges are reversed,
and the root-edge is flipped to the other side of the outer face
(in fact it is more convenient to forget about the root-edge here). Such
bipolar orientations are called \emph{pole-symmetric}.  The whole
discussion on symmetric structures is summarized in the following
proposition and illustrated in Figure~\ref{fig:symmetriesFig}.

\begin{proposition}
The number $\thetasym$ counts
\begin{itemize}
\item
pole-symmetric 2-orientations with $k+1$ white vertices and $\ell +1$ black vertices,
\item
pole-symmetric separating decompositions and 2-book embeddings with $k+1$
            white vertices and $\ell+1$ black vertices,
\item
symmetric twin pairs of full binary trees with $k+1$ left leaves and $\ell+1$ right leaves,
\item
pole-symmetric rectangulations of $X_n$ with $k$ horizontal and $\ell$ vertical segments
\item
symmetric Baxter permutations of $k+\ell+1$ with $k$ descents and $\ell$ rises,
\item
pole-symmetric plane bipolar orientations with $k$ inner faces and $\ell$ non-pole vertices.
\end{itemize}
\end{proposition}

% *************************************************************
\section{Schnyder Families}\label{sec:schnyder}
% *************************************************************

Baxter numbers count 2-orientations on quadrangulations and several
other structures. We now turn to a family of structures which are
equinumerous with 3-orientations of plane triangulations.

Consider a plane triangulation $T$, i.e., a maximal plane
graph, with $n$ vertices and three special vertices $a_1,a_2,a_3$
in clockwise order around the outer face.

\begin{definition}
An orientation of the inner edges of $T$ is a \term{3-orientation} if
every inner vertex has outdegree three. From the count of edges it follows
that the special vertices $a_i$ are sinks in every 3-orientation.
\end{definition}

\begin{definition}
An orientation and coloring of the inner edges of $T$ with colors red, green and
blue is a \term{Schnyder wood} if:
\Item(1) All edges incident to $a_1$ are red, all edges incident to $a_2$ are green
and all edges incident to $a_3$ are blue.
\Item(2) Every inner vertex $v$ has three outgoing edges colored red, green
and blue in clockwise order. All the incoming edges in an interval between
two outgoing edges are colored with the third color, see Figure~\ref{fig:vertex-3-cond}.
\end{definition}
\vskip-8mm
%%%%%%%%%%%%%%%%%%%%%%%%%%%%%%%%%%%%%%%%%%%%%%%%%%%
%%
% in einem figure environment mit caption
   \calc_figscale{30}
    \begin{figure}[htb]
    \centerline{\input{\path/vertex-3-cond.pstex_t}}
    \caption{\label{fig:vertex-3-cond}}
    \end{figure}
    VC
{Schnyder's edge coloring rule.}
%%
%%%%%%%%%%%%%%%%%%%%%%%%%%%%%%%%%%%%%%%%%%%%%%%%%%%
\ni

\begin{theorem}\label{thm:sw+3or}
Let $T$ be a plane triangulation with outer vertices $a_1,a_2,a_3$.
Schnyder woods and 3-orientations of $T$ are in
bijection.
\end{theorem}

The proof is very similar to the proof of Theorem~\ref{thm:sd+2or}.
Given an edge $e$ which is incoming at~$v$, we can classify the outgoing edges
at~$v$ as left, straight and right. Define the straight-path of an edge
as the path which always takes the straight outgoing edge. A count and Euler's formula
shows that every straight-path ends in a special vertex.
The special vertex where a straight-path ends
determines the color of all the edges along the path.
It can also be shown that two
straight-paths starting at a vertex do not rejoin. This implies that
the coloring of the orientation is a Schnyder wood.

From this proof it follows that the local properties
(1) and (2) of Schnyder woods imply:
\Item(3) \textit{The edges of each color form a tree rooted
at a special vertex and spanning all the inner vertices.}

\ni
Recall that in the case of separating decompositions we also
found the tree decomposition being implied by local conditions
(c.f. item (3) after the proof of Theorem~\ref{thm:sd+2or}).

\Note
Schnyder woods were introduced by Schnyder in~\cite{s-pgpd-89} and
\cite{s-epgg-90}. They have numerous applications in the context of
graph drawing, e.g.,~\cite{br-scdpg-05,bfm-cd3cpg-07,lls-icvrpgsr-04},
dimension theory for orders, graphs and polytopes,
e.g.,~\cite{s-pgpd-89,bt-odcp-93,fk-os-06}, enumeration and encoding of
planar structures, e.g.,~\cite{ps-ocst-03,fps-ccs3pg-04}.
The connection with 3-orientations was found by de Fraysseix and
Ossona de Mendez~\cite{fm-tao-01}.
\EndNote

\ni
The aim of this section is to prove the following theorem
of Bonichon.

\begin{theorem}\label{thm:count_schnyder}
The total number of Schnyder woods
on triangulations with $n+3$ vertices is
$$
V_n = C_{n+2}\,C_n - C_{n+1}^2 = \frac{6\,(2n)!\; (2n+2)!}{n!\; (n+1)!\; (n+2)!\; (n+3)!}
$$
where $C_n$ is the Catalan number.
\end{theorem}

Before going into details we outline the proof.
We first show a bijection between Schnyder woods and
a special class of bipolar orientations of plane graphs.
We trace these bipolar orientations through the bijection with
separating decompositions, twin pairs of trees and triples of
non-intersecting paths. Two of the three paths turn out to be equal
and the remaining pair is a non-crossing pair of Dyck paths.
This implies the formula.

\Note
The original proof, Bonichon~\cite{b-brmpgpdp-05}, and a more
recent simplified version, Bernardi and Bonichon~\cite{bb-cirt-06},
are also  based on a bijection between Schnyder woods and pairs of
non-crossing Dyck paths. In~\cite{bb-cirt-06} the authors also enumerate
special classes of Schnyder woods.

\FreeItem{}
Little is known about the number of Schnyder
woods of a fixed triangulation.
In~\cite{fz-nao-07} it is shown that the maximal number of Schnyder
woods a triangulation on $n$ vertices can have is asymptotically
between $2,37^n$ and $3,56^n$. As with 2-orientations, the
computational complexity of the counting problem is unkown.
\EndNote

\begin{proposition}
There is a bijection between Schnyder woods on triangulations with
$n+3$ vertices and bipolar orientations of maps with $n+2$ vertices
and the special property:
\Item($\star$) The right side of every bounded face is of length two.
\end{proposition}

\Proof
Let $T$ be a triangulation with a Schnyder wood $S$. With
$(T,S)$ we associate a pair $(M,B)$, where $M$ is a subgraph of $T$
and $B$ a bipolar orientation on $M$.  The construction is in two
steps.  First we delete the edges of the green tree in $S$ and the special
vertex of that tree, i.e., $a_2$, as well as the two outer edges incident to $a_2$,
from the graph, the resulting
graph is $M$.  Then we revert the
orientation of all blue edges and orient the edge $\{a_3,a_1\}$ from
$a_3$ to $a_1$, this is the orientation $B$. Figure~\ref{fig:schnyder-bip}
shows an example.

%%%%%%%%%%%%%%%%%%%%%%%%%%%%%%%%%%%%%%%%%%%%%%%%%%%
%%
% in einem figure environment mit caption
   \calc_figscale{20}
    \begin{figure}[htb]
    \centerline{\input{\path/schnyder-bip.pstex_t}}
    \caption{\label{fig:schnyder-bip}}
    \end{figure}
    VC
{A Schnyder wood and the corresponding bipolar orientation.}
%%
%%%%%%%%%%%%%%%%%%%%%%%%%%%%%%%%%%%%%%%%%%%%%%%%%%%

The orientation $B$ has $a_3$ as unique source and $a_1$ as unique
sink. To show that it is bipolar we verify properties V and F.
Property V requires that at a vertex $v\neq s,t$ the edges partition
into nonempty intervals of incoming and outgoing edges, this
is immediate from the edge coloring rule (2) and the construction of $B$.

For Property F consider a bounded face $f$ of $M$.  Suppose that $f$
is of degree $> 3$, then there had been some green edges triangulating
the interior of $f$.  The coloring rule for the vertices on the
boundary of $f$ implies that these green edges form a fan, as
indicated in Figure~\ref{fig:green-edges}. From the green edges we
recover, again with the coloring rule, the orientation of the boundary
edges of $f$ in $B$: the neighbors of the tip vertex of the green
edges are the unique source and sink of $f$. This also implies that
the right side of $f$ is of length two, i.e., $(\star)$.

If $f$ is a triangle, then two of its edges are of the same color, say
red. The coloring rule implies that these two edges point to their
common vertex, whence the triangle has unique source and sink. Since
the transitive vertex of $f$ has a green outgoing edge in $S$, it is on
the right side and $(\star)$ also holds for $f$.

%%%%%%%%%%%%%%%%%%%%%%%%%%%%%%%%%%%%%%%%%%%%%%%%%%%
%%
% in einem figure environment mit caption
   \calc_figscale{17}
    \begin{figure}[htb]
    \centerline{\input{\path/green-edges.pstex_t}}
    \caption{\label{fig:green-edges}}
    \end{figure}
    VC
{From a generic face in $S$ to $B$ and back.}
%%
%%%%%%%%%%%%%%%%%%%%%%%%%%%%%%%%%%%%%%%%%%%%%%%%%%%

For the converse mapping, consider a pair $(M,B)$ such that $(\star)$ holds.
Every vertex $v\neq s,t$ has a unique face where it belongs to the right
side. This allows us to identify the red and the blue outgoing edges of $v$.
Property $(\star)$ warrants that there is no conflict.
The green edges are the edges triangulating faces of larger degree
together with edges reaching from the right border to the additional outer vertex $a_2$.
This yields a unique Schnyder wood on a triangulation.
\qed

Given a plane bipolar orientation $(M,B)$ with $n+2$ vertices and the $(\star)$ property, we
apply the bijection from Proposition~\ref{thm:bi} to obtain a quadrangulation
$Q$ with a separating decomposition. Property  $(\star)$ is equivalent to
\Item($\star'$) Every white vertex (except the rightmost one) has a
unique incoming edge in the blue tree.
\smallskip

\ni
In particular it follows that
there is a matching between vertices $v\neq s,t$ and bounded faces of $M$,
hence, in $Q$ there are $n+2$ black and $n+1$ white vertices.

The separating decomposition of $Q$ yields twin-alternating trees with
$n+1$ black and $n$ white vertices (Theorem~\ref{thm:2or+ta}).  From
the twin-alternating pair we get to a twin binary pair of trees with
$n+1$ black and $n$ white vertices (Theorem~\ref{thm:ta+tb}).  This
pair of trees yields a triple of non-intersecting paths
(Theorem~\ref{thm:twins-3paths}).  Figure~\ref{fig:schnyder-to-paths}
shows an example of the sequence of transformations.

%%%%%%%%%%%%%%%%%%%%%%%%%%%%%%%%%%%%%%%%%%%%%%%%%%%
%%
% in einem figure environment mit caption
   \calc_figscale{12}
    \begin{figure}[htb]
    \centerline{\input{\path/schnyder-to-paths.pstex_t}}
    \caption{\label{fig:schnyder-to-paths}}
    \end{figure}
    VC
{From a Schnyder wood to three strings.}
%%
%%%%%%%%%%%%%%%%%%%%%%%%%%%%%%%%%%%%%%%%%%%%%%%%%%%

From ($\star'$) we get some crucial
properties of the fingerprint and the bodyprints of the
blue tree $T^b$ and the red tree $T^r$.

\Fact 1. If we add a leading 1 to the reduced fingerprint $\fX$, then
    we obtain a Dyck word; in symbols $(01)^{n} \leqdom 1+ \fX$.
\smallskip

\Proof
It is better to think of $1+\fX$ as the fingerprint $\fp^b$ of the blue tree after
removal of the last 0. Property ($\star'$) implies that
there is a matching between all 1's and all but the last 0's
in the $\fp^b$, such that each 1 is matched to a 0 further to the
right.\qedclaim

\Fact 2. The fingerprint uniquely determines the bodyprint of the
blue tree, precisely $\ovl{\bp^b} = 1+\fX$.
\smallskip

\Proof From ($\star'$) it follows that $\fp^b_i = 1$ implies
$\bp^b_{i+1}=0$. Since $\fp^b$ has $n$ entries 1 and $\bp^b$ that
same number of 0's, it follows that $\bp^b$ is determined by
$\fp^b$. \qedclaim

Let $\fp^* = 1+ \fX$ and $\bp^* = 1 + \bX^r$; then $(01)^{n}
\leqdom \fp^* \leqdom \bp^*$. We omit the proof that actually
every pair $(\fp^*,\bp^*)$ of 0,1 strings from $\binom{2n}{n}$ with
these properties comes from a unique Schnyder wood on a
triangulation with $n+3$ vertices. Translating the resulting
bijection with strings into the language of paths we obtain:

\begin{theorem}
There is a bijection between Schnyder woods on triangulations with $n+3 $ vertices and
pairs $(P_1,P_2)$ of
non-intersecting upright lattice paths, where $P_1$ is from
$(0,0)$ to $(n,n)$, $P_2$ is from $(1,-1)$ to $(n+1,n-1)$, and the paths
stay weakly below the diagonal, i.e., they avoid all points $(x,y)$ with
$y > x$.
\end{theorem}

For the actual counting of Schnyder woods we again apply the lemma of Gessel and Viennot.
The entry $A_{i,j}$ in the matrix is the number of paths from
the start of $P_i$ to the end of $P_j$ staying weakly below the diagonal.
The reflection principle of D.~Andr\'e allows us to write these numbers as differences of
binomials.

\begin{proposition}
\def\centbin#1{{2n \choose #1}}
  The number of Schnyder woods on triangulations with $n+3$ vertices is
   $$
    \det\begin{pmatrix}
        \centbin{n} - \centbin{n-1} && \centbin{n+1} - \centbin{n-2} \\[1mm]
        \centbin{n-1} - \centbin{n-2} && \centbin{n} - \centbin{n-3}
         \end{pmatrix}
= \frac{6\,(2n)!\; (2n+2)!}{n!\; (n+1)!\; (n+2)!\; (n+3)!}
$$
\end{proposition}

\noindent\emph{Acknowledgements.} Mireille Bousquet-M\'elou and
Nicolas Bonichon are greatly thanked for fruitful discussions.

% *************************************************************
% *************************************************************

\bibliography{barc}
\bibliographystyle{my-siam}

% *************************************************************
% *************************************************************

% *************************************************************
\end{document}